# LE LETTERE DI
# ENRICO BETTI
## (DAL 1860 AL 1886)

# CONSERVATE PRESSO
# L'ISTITUTO MAZZINIANO DI GENOVA

a cura di Paola Testi Saltini

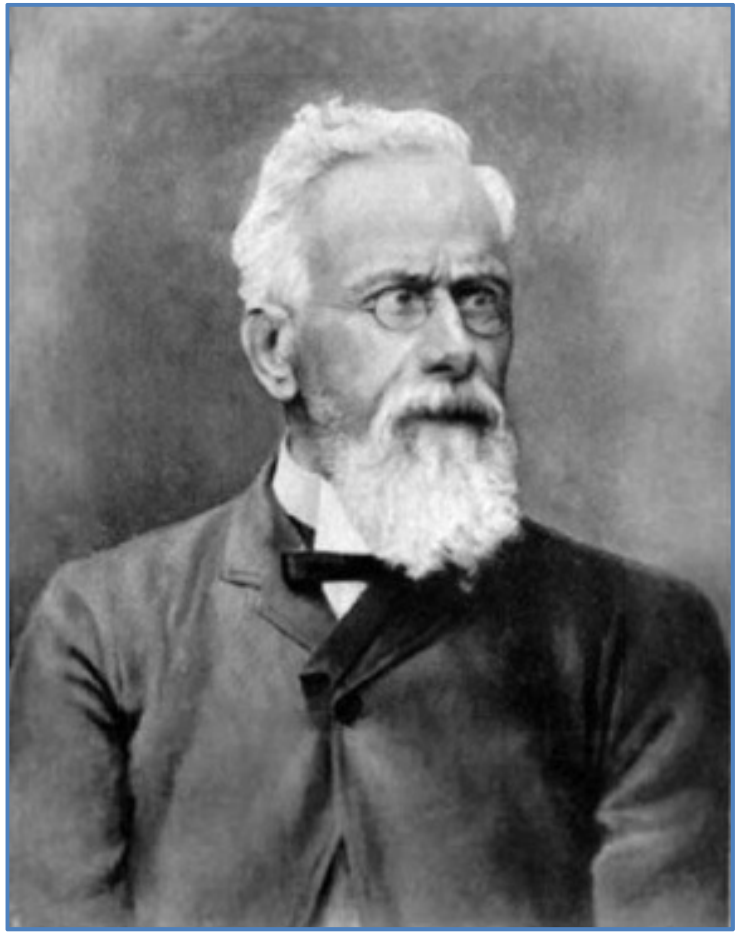

# Indice

In copertina: Enrico Betti (https://commons.wikimedia.org/wiki/File:Betti_Enrico.jpg) e sua firma autografa.

## Presentazione della corrispondenza

La corrispondenza qui riprodotta è composta da tutte le lettere spedite o ricevute da Enrico Betti e conservate presso l'Archivio Mazziniano di Genova così come sono state fino a oggi reperite. Si tratta di: 79 lettere scritte da Betti a Luigi Cremona tra il 1860 e il 1886; due lettere scritte da Betti a Eugenio Beltrami, una a Pietro Blaserna e una invece indirizzata a Betti da uno dei figli di Gaetano Giorgini che risulta però non meglio identificato.

Le lettere a Cremona trattano di questioni legate alle istituzioni ai cui lavori entrambi partecipavano, come l'Accademia dei XL o il Ministero di Pubblica Istruzione e di questioni legate alla vita universitaria. I due studiosi trattano spesso di scritti in corso di pubblicazione sugli *Atti dell'Accademia dei Lincei* oppure sul *Giornale di Matematiche*. Molte lettere sono dedicate a scambi di informazioni riguardo al volume E. Betti, F. Brioschi, *Gli elementi di Euclide con note aggiunte ed esercizi ad uso de' Ginnasi e de' licei* (Firenze, Le Monnier, 1867) al quale Cremona diede un grande apporto pur senza figurare tra gli autori.

Il rapporto che traspare dalla corrispondenza è quello tra due studiosi ad alto livello nel panorama scientifico nazionale che si stimano molto l'un l'altro e che si confrontano su parecchie questioni, sia politiche che scientifiche cercando appoggio e riferimento. Cremona è di qualche anno più giovane di Betti e, dalla prima all'ultima lettera, si avverte da parte sua il rispetto per il più anziano. È vero che Cremona diventerà Senatore prima di Betti, ma comunque l'esperienza di quest'ultimo, soprattutto per quanto riguarda le situazioni in ambito politico, viene sempre riconosciuta da Cremona.

È possibile ricostruire quasi integralmente la corrispondenza tra i due scienziati se si fa riferimento al volume [Menghini, 1996] che contiene le lettere di Cremona a Betti conservate nell'Archivio Betti della Svuola Normale Superiore di Pisa, al quale rimando anche per un'introduzione più ampia e dettagliata al carteggio e che ho usato come base di questo lavoro.

Nella pagina che segue ho inserito il breve annuncio necrologico di Betti preparato da Francesco Brioschi per gli *Annali di Matematica* in modo da consentire al lettore un orientamento spazio-temporale nella vita di Betti.

## Criteri di edizione

Poiché le lettere sono conservate in buono stato e la grafia è abbastanza chiara, la trascrizione è conforme al testo originale, sia nell'uso dei rientri che della punteggiatura, degli apostrofi e degli accenti. In particolare spesso manca l'apostrofo tra l'articolo indeterminativo e il sostantivo di genere femminile.
Si è mantenuto anche il tipo di sottolineatura del manoscritto.
Si sono mantenuti termini oggi desueti, ma in uso all'epoca, sono stati segnalati invece con *[sic!]* termini e periodi che appaiono scorretti o comunque di senso poco chiaro. I nomi propri sono stati riprodotti come scritti dall'autore, anche se inesatti e sono stati corretti soltanto in nota.
Non è stato possibile decifrare alcune parole e in tali casi si è lasciata tra parentesi quadre un'interpretazione plausibile, oppure si è sostituita la parola con il segno *[...]*.

Tutte le scritte in corsivo tra parentesi quadre sono state inserite dal curatore.

La lettera n. 48 riporta l'anno 1872, ma è quasi certamente del 1878 e, in tal caso, andrebbe a inserirsi dopo la n. 67; la n. 78 potrebbe essere del 1874 e si inserirebbe dopo la n. 56; la lettera n. 79 potrebbe essere del dicembre 1875 e quindi inserirsi dopo la n. 57.

## Enrico Betti (Pistoia, 21/10/1823 - Soiana, Pisa, 11/8/1892)

ANNUNCIO NECROLOGICO
(*Annali di Matematica*, v. XX, 1892, p. 256)

Il giorno 16 dello scorso mese di Agosto furono rese in Pisa solenni onoranze funebri ad uno degli uomini che per le sue eminenti qualità intellettuali e morali maggiormente onorarono l'Italia nella seconda metà di questo secolo, ad ENRICO BETTI.

Conobbi personalmente il Betti nelle ferie Pasquali dell'anno 1858, allorquando convenuti in Genova coll'ottimo GENOCCHI nella casa ospitale di PLACIDO TARDY, iniziammo le trattative col TORTOLINI per la pubblicazione di questi *Annali*.

Già il TORTOLINI con giusto presentimento del risveglio in quel momento degli studi matematici in Italia aveva pubblicato otto volumi dei suoi *Annali di matematica e di fisica*; ed era stato grandemente ajutato nella impresa da geometri italiani e specialmente dai su nominati.

In questi volumi si trovano quelle memorie del BETTI *Sulla risoluzione delle equazioni algebriche* (Vol. 3.°), *Sopra l'abbassamento delle equazioni modulari delle funzioni ellittiche* (Vol. 4.°), *Sopra la teorica delle sostituzioni* (Vol. 6.°), nelle quali egli, pel primo, penetrando negli alti ed astrusi concepimenti di GALOIS sulla teoria delle equazioni algebriche, pose in luce la esattezza delle asserzioni contenute nella celebre lettera a CHÉVALIER, iniziando così il movimento in quelle ricerche alle quali oggi ancora è dedicata l'attività di eminenti geometri.

Divenuto collaboratore principale di questi *Annali* il BETTI vi pubblicò dapprima varie memorie relative alla teorica delle forme, ed una monografia sulle funzioni ellittiche (Anno 1860), nella quale l'importante teoria è esposta con originalità di metodo.

Dal 1862 in poi la sua attività scientifica può dirsi tutta rivolta agli ardui problemi della fisica matematica, ed anche in questo campo la sua acuta intelligenza, le profonde sue cognizioni nell'analisi matematica, si rivelarono ben presto nei risultati da lui ottenuti. Oltre venti memorie sue sopra questa disciplina sono contenute in questi *Annali*, nel *Nuovo Cimento*, nelle *Memorie* della Società Italiana dei Quaranta, nei *Rendiconti* della R. Accademia dei Lincei. Infine a questo periodo della sua vita devesi quel suo libro sul *Potenziale*, tradotto in Germania e tanto apprezzato dai cultori di questi studi.

Il BETTI non ebbe lunga esistenza, essendo nato nei pressi di Pistoja, il 21 ottobre 1823. Egli fu insegnante esimio e caldo amico dei suoi discepoli. A lui principalmente spetta l'onore di avere, colla sua parola e col suo esempio, fondata quella Scuola matematica dell'Università di Pisa dalla quale sono usciti in molta parte i giovani valorosi che hanno mantenuto e mantengono alto il nome italiano negli studi matematici. È a questi giovani geometri che io, rammentando la cospicua partecipazione del BETTI nell'iniziare la pubblicazione di questi *Annali* e nel collaborarvi, raccomando il loro avvenire.

Ottobre 1892.

F. BRIOSCHI

Per l'elenco dettagliato dei lavori scientifici di Betti si veda il volume I delle *Opere di Enrico Betti*, Hoepli, Milano, 1903, pp. V-VII, disponibile anche in rete.

# Il carteggio

### *Lettere a Luigi Cremona*

**1**

Stimatissimo Sig<sup>r</sup> Professore,

Ho ricevuto dal comune amico Casorati il suo ritratto che ho avuto gratissimo, e la ringrazio. La ringrazio anche del suo pregevole opuscolo sopra la Geometria che ebbe la gentilezza di inviarmi un mese fa.

Mi affretto a soddisfare il suo desiderio di avere il mio ritratto, che mi ha espresso il Casorati, e lo fo molto volentieri perché mi è veramente grato di potere stringere vincoli di amicizia con Lei che stimo molto per il suo amore alla Scienza, per l'abilità che in essa dimostra con i lavori che va pubblicando, e per il suo carattere che merita l'amicizia del mio ottimo amico Brioschi.

Non so quando potrò fare una corsa a Bologna, ma non è probabile prima dell'estate. Nonostante sarà un motivo forte per indurmi a ciò il desiderio di conoscerla personalmente.

Mi confermo pieno di stima.

Suo aff<sup>mo</sup> Amico
Enrico Betti

Pisa, 19 Novembre 1860

**2**

Torino, 15 Febbrajo 1861

Amico carissimo,

Non ho risposto prima alla vostra lettera, perché quando è arrivata a Pisa io ero già partito ed è stata respinta a Firenze e di lì a Torino, ed ha sofferto così un ritardo di parecchi giorni prima di arrivarmi.

Io aveva ricevuto le tre Copie della vostra prolusione,[1] due delle quali detti a Mossotti e Novi come accennava che io dovessi fare l'indirizzo.

La Cattedra di Geometria Superiore nell'Università di Pisa fu riunita nominalmente in quella di Analisi superiore. Ridolfi[2] che aveva in corso un numero stragrande di Professori per Firenze, non volle aumentare di questa Cattedra l'Università di Pisa nonostante che qualcuno ne mostrasse la convenienza. Vedremo come si potrà ideare il nuovo ordinamento; è probabile che risulti la necessità di questa Cattedra a Pisa o a Firenze. Allora io mi auguro di tutto cuore di potervi avere per Collega.

Io sono qua fino dal 4 Febbrajo. La Commissione ha fatto soltanto tre adunanze, ed io ora son costretto a starmene qui senza studiare e facendo niente; ciò che mi noja moltissimo. Speriamo che arrivati quei Commissari che mancano tuttavia, si procederà con attività e si concluderà qualche cosa.

Salutatemi Chelini e ditegli che ho ricevuto oggi la sua Memoria e lo ringrazio. Unita ve n'era una per Mossotti, che subito vo a respingere a Pisa.

Novi mi disse che desideravi *[sic!]* avere un esemplare della mia Monografia sopra le funzioni ellittiche.[3] Appena riceverò gli esemplari che ho fatto tirare a parte, ve ne manderò subito uno. Resta a pubblicarsi la seconda parte che esposi l'anno passato nelle lezioni date nell'Università, ma non l'ho ancora ordinata per la stampa, e sebbene sia cosa che non richiederà molto tempo, pure finché sto qua a Torino difficilmente potrò fare qualche cosa per piccola che sia.

Credetemi sempre

Vostro aff<sup>mo</sup> Amico
Enrico Betti

---

[1] Si tratta della "Prolusione ad un corso di geometria superiore, letta nell'Università di Bologna. Novembre, 1860", *Il Politecnico*, 1861, pp. 22-42.
[2] Potrebbe trattarsi del senatore Cosimo Ridolfi.
[3] E. Betti, "La teoria delle funzioni ellittiche", *Annali di Matematica pura ed applicata*, 1859-1860.

**3**
*Carta intestata:* Camera dei Deputati

Torino, 12 Marzo 1863

Carissimo amico,

Appena ricevuta la tua lettera ho parlato del tuo affare col ministro Amari.[4] Egli mi ha detto che ora occupato com'è per la diffusione del Bilancio della istruzione alla Camera non può attendere, ma che ne riparleremo con più comodo dopo che sarà terminata questa discussione.

Io tornerò presto in Toscana, e là spero che ci vedremo al tuo ritorno. Vedremo anche meglio che cosa potremo proporre al Ministro per rimediare. Probabilmente ci sarà la difficoltà economica, perché alla Camera si fanno molte resecazioni sul Bilancio, e il Ministro forse sarà obbligato a far dei risparmi sul personale. Nonostante non ti salti in testa l'idea di dare la dimissione [sic!]. L'idea che un insegnamento superiore vi debba essere sopra ogni parte dello scibile in qualche Università o Istituto superiore è generale, e tu che hai tanto valore nella Geometria superiore non puoi dubitare che prima o dopo non possa aver tempo di dare questo insegnamento. Io prima di partire di quà ne riparlerò col ministro, e poi ci vedremo a Pisa o a Firenze. A Pisa io abito in casa [Fiorelli] dietro l'Arcivescovado 2°. Piano. Se per avventura tu arrivassi a Pisa ed io non vi fossi, il che potrebbe essere se salterò le vacanze di Pasqua, sarei a Firenze in Via dell'Oriuolo #° 37. 1° Piano. Salutami Battaglini, Padula e De Gasparis ed ama

il tuo amico aff$^{mo}$
Enrico Betti

**4**

Pisa, 9 Aprile 1863

Amico carissimo,

Scrissi subito al Ministro,[5] ed Egli mi ha risposto immediatamente.
Ecco le sue parole = Farò il mutamento che voi credete conveniente nella posizione del Prof$^e$ Cremona, il quale me ne ha già scritto ed io voglio onorarlo e soddisfarlo in ogni modo come Egli merita = Mi dice poi che fino ad ora era rimasto incerto soltanto per l'autorità scientifica dei due che avevano fatto il Decreto precedente a tuo riguardo. Io sono lietissimo di aver contribuito in parte a che ti sia fatta giustizia, siano soddisfatti i tuoi desideri, e ti sia lasciato maggiore agio agli studi della Scienza che coltivi con passione ed amore.

È uscito alla luce il Ricordo del Prof. Mossotti scritto dal Bicchierai.[6] Si venda al prezzo di una Lira a profitti dei Monumenti Mossotti e Salvagnoli.[7] Te ne devo mandare un certo numero di copie? Quante credi che se ne possa spacciare a Bologna?

Salutami i comuni amici ed ama

il tuo Amico
Enrico Betti

---

[4] Si vedano [Menghini, 1996], lettere nn. 8-11, pp 13-14. Cremona, che si trovava a Cosenza per un giro d'ispezione, era stato spostato dalla cattedra di Geometria superiore a quella di Geometria analitica e descrittiva presso l'Università di Bologna senza venire consultato dal Ministro. Sdegnato, il geometra aveva minacciato le dimissioni nella lettera a Betti del 4 marzo. Già il 5 aprile gli scriveva in altro tono e l'11 aprile riferisce al collega, ringraziandolo per l'intervento, di aver ricevuto dal Ministro la notizia del reintegro sulla cattedra di insegnamento superiore unita a "cose assai lusinghiere" e all'invito a recarsi in Sicilia per un altro giro d'ispezione delle scuole secondarie.
[5] Michele Amari.
[6] Z. Bicchierai, "Ricordo del prof. O.M. Mossotti", *Gazzetta di Firenze*, n. 60.
[7] Il 27 aprile 1873 nel camposanto di Pisa venne dedicata una statua a Vincenzo Salvagnoli.

**5**

<div style="text-align: right;">Pisa, 15 Aprile 1863</div>

Amico Carissimo,

Riceverai 20 Copie del Ricordo del Mossotti. Ho creduto più sicuro indirizzarle a te, che le farai recapitare al librajo che presta senza interesse. Quando partirai per la Sicilia[8] lascia pure incaricato il Capellini: però non si potrà mettere nel numero dei *[?]* componenti la Giunta promotrice che si è annunziata definitivamente costituita. Credo però che il Capellini che conosceva personalmente il Mossotti si presterà volentieri in ogni modo. Anche Cecco Magni che conosceva bene il Mossotti credo che si presterà volentieri a raccogliere firme e denaro. Quando sarà raccolto un certo numero di firme il Capellini me ne manderà nota perché si possano pubblicare nel Giornale di Pisa secondo la promessa.

Ricevi i saluti di Novi e di Lavagna, scrivimi prima di partire e dimmi la strada che fai perché se ci potessimo vedere sarebbe molto grato.

Credimi sempre

<div style="text-align: right;">il tuo Aff. Amico<br>Enrico Betti</div>

**6**

<div style="text-align: right;">Spezia, 24 Agosto 1863</div>

Amico carissimo,

Il Beltrami prima di avere ricevuto la tua risposta mi ha scritto dicendomi che non accetta. Tra le ragioni che adduce una delle più forti è che Egli crede molto più meritevole Delfino Codazzi. Ma di questo Codazzi io ho informazioni che me lo farebbero credere un pessimo collega, e io non potrei pensare a proporlo. Che cosa debbo fare? passare sopra agli scrupoli del Beltrami, oppure insistere poiché la Cattedra sia messa al concorso?
Io credo che potrei proporre al Ministro qualcuno che potrebbe far buona riuscita nella Cattedra che lascerebbe il Beltrami.

Io dimani andrò a Pisa; tra due o tre giorni sarò di ritorno qui, ove mi tratterrò tutto Agosto e forse qualche giorno di Settembre. Avrò molto piacere se potrò vederti alla Spezia.

Ama il

<div style="text-align: right;">tuo amico<br>Enrico Betti</div>

**7**

<div style="text-align: right;">Spezia, 29 Agosto 1863</div>

Amico carissimo

Ho scritto al Ministro[9] proponendo Beltrami, di cui già gli aveva parlato.
Bisognerebbe che Beltrami, appena nominato, venisse a Firenze. Là all'Osservatorio con Donati prenderebbe le cognizioni pratiche necessarie, ed io gli mostrerei il corso che faceva Mossotti e qualche consiglio potrei dargli quanto al modo di regolarsi. Quando avrò avviso dal Ministro che non vi sono difficoltà scriverò io stesso al Beltrami. Intanto se hai occasione di scrivergli potresti dargli questo consiglio.

---

[8] Cremona, su incarico del Ministero, stava iniziando un "giro d'ispezione per le Puglie e per le Calabrie" come scrive in [Menghini, 1996], lettera n. 8, p. 13. Si vedano anche [Menghini, 1996], lettere nn. 12 e 13, p. 15.
[9] Michele Amari.

Dubito di vederti alla Spezia. Io Lunedì mattina parto. Starò qualche giorno a Pisa, poi andrò a Firenze dove mi tratterrò il Settembre e l'Ottobre. Per ora se hai occasione di scrivermi, indirizzami la lettera a Pisa.
 Ama il

<div style="text-align:right">tuo amico<br>E Betti</div>

**8**
*Carta intestata:* Camera dei Deputati

<div style="text-align:right">Torino, 4 Dicembre 1863</div>

Mio caro amico,

 Mi ha detto il Ministro[10] che il Padula farà lezione da sé; quindi per ora non vi è luogo dar quella Cattedra. Quanto ai libri di testo puoi mandare il tuo lavoro senza intenderti prima cogli altri due Commissari.
 Quanto al progetto della sotto commissione per la riforma della istruzione superiore mi disse il Piria che me ne avrebbe mandata una copia ma non l'ho ancora veduta. L'ho chiesta al Bonghi il quale mi ha detto che non ne ha. Potrò averla dal Ministro.
 Questo progetto però è nato morto, perché gli è contrario lo stesso Ministro. Darebbe alle Università quello che hanno presentemente e le lascerebbe governarsi da sé il tutto e per tutto.
 Ebbi molto piacere che la lettera del Ministro inducesse il Chelini a desistere dall'idea di dimettersi o lo tranquillizzasse. Salutalo da parte mia.
 Io probabilmente partirò dimani sera di qui per Toscana. Ma per quanto grande sia il desiderio di vederti e parlare un poco con te non posso dirti che tu venga alla stazione a trovarmi, perché arriverò costì la mattina alle 5 e 20 minuti e ripartirò alle 6 e 21, e a questo freddo vi è di prendere un malanno. Sicché ci vedremo quando ripasserò per tornare a Torino. Ora che è aperta la strada fino a Pracchia,[11] ho abbandonato la via di mare, e passerò spesso da Bologna.
 Salutami i comuni amici, scrivimi a Pisa e ama

<div style="text-align:right">il tuo Amico<br>Enrico Betti</div>

**9**
*Carta intestata:* Camera dei Deputati

<div style="text-align:right">Torino, 7 Marzo 1864</div>

Mio caro amico,

 Non ho potuto vedere Amari prima di oggi; ma ho da darti eccellenti nuove.[12] Il Ministro non aveva niente a dirmi a te sfavorevole, anzi è pienamente nelle tue vedute. Mi ha detto che Egli ha deliberato che Fiorini debba riprendere le lezioni, che stasera sarà portato l'affare in Consiglio dei Ministri, e che se il Consiglio dei Ministri la dà vinta agli scolari e al Rettore,[13] egli lascia il portafoglio. Le cose stanno come ti dissi alla stazione.
 Dirai a Magni che gli ho parlato della Cattedra di Pisa e che mi ha detto che non ha intenzione di conferirla a nessuno.
 Salutami i comuni amici e riama

<div style="text-align:right">il tuo Amico E Betti</div>

---

[10] Michele Amari.
[11] Comune in provincia di Pistoia.
[12] Si vedano [Menghini, 1996], lettere nn. 24-27, pp. 21-23.
[13] Antonio Montanari.

**10**
*Carta intestata:* Camera dei Deputati

Torino, 18 Luglio 1864

Amico carissimo,

Saprai forse che Padula invitato dal Ministro[14] ad ottare tra gl'impieghi che occupava, ha scelto la Cattedra della Università e quella del Collegio di Marina, ed ha rinunziato alla Direzione e alla Cattedra della Scuola di Applicazione. Io ne ho parlato col Ministro il quale è stato dispiacente che questa scelta del Padula non lasci vacante la Cattedra Universitaria, e quindi renda impossibile per ora la combinazione che Egli aveva stabilito per il tuo traslocamento a Napoli. Rimarrebbe a credere se fosse possibile un'altra combinazione; ma il Ministro ha intenzione di lasciare il portafoglio appena gli sia possibile e non vuole affrettarsi a far nomine prima del tempo quasi che volesse *[?]* il compito del suo successore, che non saprei dirti chi potrà essere e quando accadrà questa mutazione.

Stasera tornerò in Toscana, ma passerò per Bologna. Quindi il ritratto di Clebsh *[sic!]*[15] mandamelo a Pisa.
Salutami Tardy. Scrivimi ed ama

il tuo Amico
Enrico Betti

**11**

Pisa, 19 Maggio 1865.

Mio caro amico,

Come ti avrà detto Magni, Mercoledì aveva fissato col Ministro di andare da Lui per parlare del progetto di cui ti scrisse Beltrami, ma avendo Egli dovuto anticipare la sua partenza per Torino non poté parlarci come avrei desiderato. Però credo che presto tornerà qua, e allora farò una gita a Firenze e lo vedrò. Intanto ti dirò che questo progetto è anche di Matteucci, e che io gli parlai espressamente di te quando lo vidi a Torino, e gli mostrai il vantaggio che se ne sarebbe ricavato trasportando l'insegnamento della Geometria Superiore dato da Te, da Bologna dove non vi è altri insegnamenti di complemento, a Firenze dove si farebbe un centro dove si coltiverebbero e insegnerebbero anche le altre parti della scienza a giovani che avrebbero per iscopo non di fare la professione d'ingegnere, ma di divenire cultori e docenti di Matematica, e come anche in questo ruolo saresti tolto dal trovarti insieme col Rettore Montanari, ed egli ne andava perfettamente d'accordo. Ora che è là il Ministro, Matteucci tratterà con Lui del Museo di Firenze e se credi anche di scrivergli da tè il tuo desiderio non farà male.

La medesima idea la ebbe anche Brioschi e ne parlò al Ministro un mese fa.

Se non riuscirà di fermare questa lezione di Matematica al Museo di Firenze, io parlerò al Ministro perché t'incarichi di venire ad insegnare a Pisa la Geometria Superiore.

Ho preso dal Novi la tua Memoria di Sylvester. Quando l'avrò letta te la manderò; se ti occorre prima scrivimelo.

Ama sempre

il tuo Amico
E Betti

---

[14] Michele Amari.
[15] Rudolf Friedrich Alfred Clebsch.

**12**

Di campagna presso Pistoia 14 Sett 1865

Mio caro Cremona,

Mi aveva detto il Brioschi che tu non desideravi di muoverti da Bologna per ora. Vi sarebbe stata qualche difficoltà perché il Ministro[16] è contrario a incaricare uno che è Professore in un luogo d'insegnare in un altro, come sarebbe stato necessario di fare. Ma forse si sarebbe superata ponendo avanti l'interesse della scuola Normale. Speriamo che dopo questo anno scolastico potremo avere occasione di trovarci in un medesimo istituto.

In una settimana del Luglio e in una dell'Agosto ci siamo riuniti a Firenze Giorgini, Bonghi, Brioschi, Villari ed io per compilare un progetto di legge sulla pubblica istruzione del quale ci aveva incaricati il Ministro. Quando ci dividemmo lasciammo l'incarico a Giorgini di formulare le poche proposte per l'insegnamento superiore e a Villari quelle per l'insegnamento primario e secondario. Ieri fui a Firenze e vidi il Bianchi[17] il quale mi disse che desiderava che queste proposte di legge gli fossero rimesse sollecitamente, perché il Ministro vuol presentarle al Parlamento. Quindi può ritenersi che sin allora terrà il portafogli Natoli e rimarrà il Bianchi Segretario generale.

Quanto alle nostre proposte non ne avrai ricavato la natura dai primi articoli che ha pubblicato il Bonghi nella Nazione, ma in gran parte dal 5° e da quelli che ha pubblicato Brioschi nella Perseveranza. Le Università attuali si riducono quasi tutte ad Atenei con venti Professori dove si fanno due anni di corso scientifico, preparatorio agli Istituti e alle Scuole Speciali: delle quali 3 per gl'Ingegneri 6 o 7 di Medicina dove sono grandi Ospedali, e di diritto almeno una per ogni Università primaria. Per formare i cultori della Scienza due o tre Istituti che darebbero il grado dottorale e di abilitazione a fare i dottori docenti.

Salutami il nostro Magni ed ama sempre

il tuo Amico
Enrico Betti

P.S. Hai notizie di Beltrami?

**13**

*Carta intestata:* R. Scuola Normale Superiore

Pisa 8 Aprile 1866

Amico carissimo,

Sabato sera ci adunammo col Ministro[18] e il Segretario Generale[19] per i rapporti da farsi sopra le Scienze Fisiche e Matematiche nell'occasione della Esibizione di Parigi. Eravamo presenti soltanto Menabrea Donati ed io. Matteucci vorrebbe che non se ne facesse nulla, perché dicono che le altre nazioni non hanno accettato l'invito della Francia. Mandò a dire che la sera non può uscire di casa. Io dimandai al Ministro se era vero che le altre nazioni non avessero accettato ed Egli mi disse che se ne sarebbe informato. Intanto disse noi facciamo, poi se non converrà mandar le relazioni a Parigi le pubblicheremo in Italia a spese del Ministero. Si vede bene che Egli si è troppo *[avanzato]*, accettando subito l'invito. Fu stabilito che il Ministero avrebbe procurato di avere notizie per tutte le pubblicazioni fatte in questi ultimi 20 anni nei diversi rami di Scienza, e che quando a noi occorressero notizie, libri ecc ci volgessimo al Ministero che ce le avrebbe procurate. Fu parlato anche del modo di scrivere le relazioni. Brioschi, tu ed io bisognerà poi che ci vediamo per fissar meglio il da farsi. Ieri sera non invitarono ne *[sic!]* Brioschi ne *[sic!]* te ne *[sic!]* Cannizzaro per non incomodarvi. La legge degli organici delle Università sarà presentata al

---

[16] Giuseppe Natoli.
[17] Nicomede Bianchi, senatore.
[18] Domenico Berti.
[19] Federico Napoli fu Segretario generale del Ministero della Pubblica Istruzione dal 6 gennaio 1866 al 15 maggio 1869.

riaprirsi del Parlamento. Allora te ne scriverò. Dini credendo che io ti avrei veduto a Firenze mi aveva dato questo ritratto da consegnarti.
    Scrivimi ed ama

il tuo Amico Betti

**14**
*Carta intestata:* R. Scuola Normale Superiore

Pisa 27 Giugno 1866

Mio caro Cremona,

    Parlammo a Firenze con Brioschi della nomina del Presidente della Società Italiana delle Scienze. Bisognerebbe scegliere uno che stesse a Firenze, e che fosse attivo e disposto a dar nuova vita alla Società. Nella Capitale del Regno non vi è Accademia delle Scienze: dovrebbe la Società dei Quaranta prendere questo posto. Ora la sede della Società è nel luogo dove risiede il Presidente: perciò Brioschi ed io ponemmo gli occhi sopra Matteucci. Che cosa ne dici tu? Avresti qualche altra persona a proporci? Scrivimi la tua opinione perché possiamo d'accordo fare un poco di propaganda, per il meglio della Società.
    Montanari parlò con Brioschi mostrandosi deciso a portarci via il nostro ottimo e bravo Beltrami. Se non ci fosse la forte ragione della salute della Madre, io contratterei con tutte le mie forze, ma di contro a quella ragione io se potrò favorirò invece di contrariare.
    Scrivimi presto ed ama

il tuo amico aff
Enrico Betti

**15**
*Carta intestata:* Camera dei Deputati

Firenze, 17 Agosto 1866

Mio caro Cremona,

    Parlai col Segretario del noto affare. Mi disse che nel caso venisse fatta una proposta dal Rettore[20] della Università la cosa sarebbe fatta subito senza difficoltà. Ma fin ora non è arrivata al Ministero nessuna proposta del Rettore. Potrebbe essere che il Rettore aspetti a farla insieme con altre.
    Se il Rettore però non volesse fare la proposta di non confermare codesto straordinario allora la cosa è meno facile. Ma la via che il Segretario[21] mi suggeriva è questa: tu stesso di cui il Ministro[22] ha molta stima potresti scrivere una lettera o a lui direttamente o a qualcuno che gliela mostrasse, e se tu credessi potrei incaricarmene anche io, nella quale tu esponessi i danni che ne risultano all'insegnamento per la incapacità di altri. Ma se il Rettore fa la proposta, quella via è preferibile.
    Il Beltrami è sempre a Padova? Scrivimi ed ama sempre

il tuo Amico
Enrico Betti

---

[20] Antonio Montanari.
[21] Potrebbe trattarsi di Federico Napoli.
[22] Domenico Berti.

**16**
*Carta intestata:* Camera dei Deputati

Firenze, 19 Agosto 1866

Mio caro Cremona,

Nel caso che il Beltrami lasciasse la Cattedra di Pisa io avrei necessità di avere il Dini in quel posto, non tanto per compensare la perdita nella Università, quanto per la Scuola Normale dove mi mancherebbe un valido ajuto. Tra quelli giovani ve ne sono che possono fare, ed hanno bisogno di essere in contatto di chi coltivi le scienze con passione come faceva il Beltrami e come farà il Dini. D'altra parte il Dini avrà più piacere di restare a Pisa dove ha la famiglia. Si potrebbe dunque nel caso che Beltrami andasse a Padova proporre la supplenza per Boschi. Se poi il Beltrami accetta Bologna la cosa va da sé, perché resta il posto di Pisa per Dini.

Scrivimi a Pisa appena vi è qualche cosa di risoluto, ed io appena l'affare sarà rimesso al Ministero andrò dal Berti e l'appoggerò quanto mi sarà possibile.

Ama sempre

il tuo Amico
Enrico Betti

**17**
*Carta intestata:* Camera dei Deputati

Firenze, 29 Agosto 1866.

Mio caro Cremona,

Oggi ho parlato con Berti del noto affare. Egli mi ha detto che è disposto a fare il meglio per l'utile delle nostre Università. Mi ha dimandato però; che cosa faremo di questo Venturi? Io gli ho risposto che potrà rimettersi a fare il Professore di Disegno, come faceva prima; che credo che a questo il Reggente[23] ci avrà pensato. Io ne ho parlato nel senso fissato precisamente, e sono rimasto col Napoli che io gli farò un appuntino per sua memoria. Me perché la cosa andasse da sé senza difficoltà e senza fatica occorrerebbe la proposta del Reggente. Pare che anche il Gatti Capo di Divisione sia informato della cosa e sia contrario, perché l'altro giorno ne parlò col Segretario.[24]

Io dimattina vo in campagna presso Pistoia. Indirizzami le lettere a Pistoia. Però verrò spesso a Firenze, e quando occorrerà che io dica un altra parola la dirò.

Salutami Cecco Magni ed ama il

Tuo Amico aff
Enrico Betti

---

[23] Antonio Montanari.
[24] Federico Napoli.

**18**

Caloria, 15 Ottobre 1866.

Mio caro Cremona,

Non vi è alcun dubbio intorno alla possibilità di differire gli esami alla Sezione di estate.

Io non ho qui la <u>Teorica delle forze che agiscono</u>[25] etc; ma parmi che nel luogo citato sia applicata la proprietà del potenziale di superficie relativa alla discontinuità della derivata. In tal caso questa proprietà sussiste evidentemente anche quando alla funzione potenziale si aggiunge una funzione potenziale P che si mantiene continua attraversando la superficie. Del resto quando ci vedremo a Firenze avrò riscontrato il libretto e ti potrò dare *[altro perché li integri all'istruzioni]* gli schiarimenti opportuni.

Scrivimi quando saprai il giorno in cui arriverai a Firenze, ché io non mancherò di venire a trovarti.

Quel tuo scolaro che voleva venire alla Scuola Normale non ha fatto ancora la dimanda.[26] Si è forse impaurito dell'esame? Non creda però che si pretendano cose straordinarie. Più che sul quanto fanno, si cerca di conoscere come sono padroni di quel che lor è stato insegnato.

Venerdì passato fui a Firenze e vidi un momento Brioschi e Matteucci, coi quali parlammo della Società Italiana.

Salutami Cecco Magni ed ama

il tuo Amico aff
Enrico Betti

**19**
*Carta intestata:* R. Scuola Normale Superiore

Pisa 8 Gennajo 1867

Mio caro Cremona,

Non vi è nessun dubbio che ora bisogna ritornare di proposito a dar vita a un periodico di Matematica e certamente sarà d'immenso vantaggio sostituire Milano a Roma per luogo di pubblicazione, e te e il Brioschi a Tortolini. Anche il pubblicare a fascicoli sciolti come il Crelle mi par fatto bene. Però la tua proposta io l'approvo pienamente e volentierissimo. Se comincerà presto la pubblicazione per il primo fascicolo potrò dar subito un articolo bibliografico sopra l'opuscolo di <u>Neumann</u> sopra le funzioni sferiche, mostrerei come col teorema di Green si può anche dimostrare la convergenza degli sviluppi per funzioni sferiche, di qualunque funzione data sulla sfera, che pare sia sfuggita a Neumann.[27]

Quanto alle relazioni a me non hanno detto nulla Berti e Napoli quando ho avuto occasione di vederli.[28] A me pare cosa mal fatta di andar noi a mostrare il nostro misero patrimonio scientifico degli ultimi 20 anni, quando chi è molto più ricco di noi modestamente vi si rifiuta. I miei colleghi relatori di qua pensano come me.

Non ho ancora esaminato i fogli del nostro povero Novi; so però che aveva cominciato a lavorare al 2° volume e doveva aver finito alcuni Capitoli.[29] Vedrò e penserò al da farsi.

---

[25] E. Betti, "Teorica delle forze che agiscono secondo la legge di Newton e sua applicazione alla elettricità statica", *Nuovo Cimento*, s. 1, XVIII, pp. 385-402; XIX, pp. 59-75, 77-95, 149-175, 357-377; XX, pp. 19-39, 121-141, 1863-64; tali articoli vennero pubblicati insieme nel 1865 per i tipi Pieraccini di Pisa.

[26] Potrebbe trattarsi di Eugenio Bertini: si vedano [Menghini, 1996], lettere n. 38 e n. 39 p. 30.

[27] Non risulta che Betti abbia recensito lavori di Neumann. Potrebbe riferirsi alla questione trattata l'unico lavoro apparso nel primo numero degli annali: E. Betti, "Sopra le funzioni sferiche", *Annali di matematica*, s. 2, v. 1, 1867, pp. 81-87. Nello stesso anno viene pubblicato: C. Neumann, *Theorie der Bessel'schen functionen: ein analogon zur theorie der Kugelfunctionen*, Teubner, 1867.

[28] Su tale questione si vedano qui la lettera n. 13, le lettere che seguono e [Menghini, 1996], lettere nn. 39-44, pp 30-35 scritte sull'argomento da Cremona.

[29] Probabilmente si tratta della seconda parte di: G. Novi, *Trattato di Algebra superiore. Parte I Analisi Algebrica*, Firenze Le Monnier, 1863.

Si vende la sua libreria. Manderò il catalogo anche a te. Intanto quello che premerebbe per ajutare la sua famiglia sarebbe di poter vendere a prezzo discreto gli ultimi 40 Volumi del Crelle; perché così la Sua Collezione completa la prenderebbe la nostra Biblioteca al prezzo di commercio.

Mancano poche centinaja di Lire per pagare lo scultore del Monumento Mossotti. Dì a Brioschi che quando viene a Firenze non si dimentichi della promessa, e vedete se potete raccapezzare qualche cosa a Milano. Il Monumento sarà messo al posto nel Maggio.

Salutami Brioschi e Casorati ed ama

il tuo amico
Betti

**20**
*Carta intestata:* R. Scuola Normale Superiore

Pisa 25 Marzo 1867

Mio caro Cremona,

Quando il Berti ci adunò a Firenze, io gli dimandai se era vero che le altre nazioni non avessero accettato l'invito della Francia. Il Berti mi rispose che non sapeva, che si sarebbe informato, ... ma questo allora non doveva preoccuparci, perché la questione se i rapporti si dovessero o non si dovessero mandare a Parigi si sarebbe risoluta dopo. Che Egli in ogni modo desiderava che questi rapporti fossero fatti, perché se non si mandassero poi a Parigi non credendolo conveniente il Governo li avrebbe fatti pubblicare a sue spese in Italia per l'utile della Scienza. Questa risposta datami presenza di Napoli, Donati e Menabrea mi faceva supporre che quando fosse verificato che le altre nazioni come la Prussia e l'Inghilterra non avessero accettato l'invito, non avremmo mandato i rapporti neppur noi, e quindi questi rimanevano lavori da farsi a comodo e come si fosse creduto meglio. Questo fatto mi pare che ci offra modo di mostrare chiaramente che noi non abbiamo mancato, e come anche nel concetto del Ministro che ci aveva nominato[30] si dovesse ritenere che non convenisse andar noi a far mostra delle nostre miserie quando chi era più ricco di noi si asteneva dal volerne far pompa.

La miglior cosa sarebbe se Brioschi andasse presto a Firenze che a voce gli esponesse tutte queste ragioni anche a nome nostro. Se no mi pare che bisognerà scrivere un altra lettera complessiva appoggiando sopra le ragioni che ci avevano dovuto rendere persuasi che quando le altre nazioni non avessero accettato noi non avremmo mandato i rapporti, e ripetendo le ragioni che possano convincere il Ministro[31] della convenienza della cosa.

Io avrò molto piacere di pubblicare qualche cosa nel primo Numero degli Annali, e ti mando un articolo sopra le funzioni sferiche.[32] Mandami a correggere le prove di stampa.

Teza e Dini ti salutano caramente. Dini ora è molto occupato perché su due corsi: Geodesia e Algebra. Egli ha trovato vari associati per gli Annali.

Salutami Brioschi e Casorati. Scrivimi ed ama sempre

il tuo Amico
Enrico Betti

---

[30] Era Domenico Berti (si veda la lettera n. 13).
[31] Cesare Correnti.
[32] Si veda la lettera n. 19.

**21**

Pistoia, 7 Settembre 1867

Mio caro Cremona,

L'amministratore della Società <u>Successori Lemonnier</u> risponde alla mia dimanda quanto al compenso per la pubblicazione di Euclide,[33] che egli crede che non vi sia ombra di dubbio sull'applicabilità nel nostro caso di quello che fu deciso in massima dalla Commissione per i libri pedagogici, che cioè si debba dare agli Autori e commentatori dei medesimi la metà degli utili netti.

I primi giorni della settimana entrante andrò a Firenze per un giorno. Da che ti lasciai a Firenze io sono stato sempre quassù nella quiete dei monti dove passeggio, insegno Euclide al mio nipote, e studio quelle funzioni delle quali ti parlai. Aspetto presto lavoro di revisione di stampa etc intorno a Euclide.

Hai notizie di Brioschi. Scrivimi presto ed ama sempre

il tuo Amico aff[mo]
Enrico Betti

**22**

Pistoia, 24 Settembre 1867

Mio caro Cremona,

Mi scrive l'Amministratore della Società <u>Lemonnier</u> che bisogna fargli sapere subito il numero delle Copie da tirarsi. Fa osservare che se non si può contare sopra lo smercio di un certo numero di copie, rimarrà poco o nulla per la metà degli utili che aumenterà in una proporzione molto rapida coll'aumento del numero delle copie vendute. Quindi o bisogna contare sopra qualche migliaio di copie in un solo anno o tener fermo l'Euclide per più. Scrivimi subito la tua opinione in proposito. Ti mando le prove di stampa corrette perché tu gli dia un ultima revisione; tra poco ti manderò le altre fino alla fine del 1° Libro. Trovo che <u>Hankel</u>[34] sostiene che in tutti i manoscritti si trovano gli assiomi 10 e 11 come postulati 4° e 5°, e trova inconcepibile come si possano porre come assiomi le proposizioni 10, 11 e 12 che sono i postulati 4, 5 e 6. Hoüel[35] sostiene il contrario. È chiaro che sono distinti dai primi 7 assiomi che valgono per tutte le grandezze non per le sole geometriche, e che sono vere nozioni comuni, perché veramente ogni uomo le ha. mentre le ultime tre sono speciali a certi enti geometrici, e una osservazione e una esperienza scientifica può solo almeno per l'assioma 11, trarle dalla esperienza. I tre postulati 1, 2, 3, sono le tre sole a cui si riducono tutte le costruzioni della geometria elementare; i tre 10, 11 e 12 sono le sole proposizioni geometriche dalle quali si deducono tutte le altre per mezzo delle nozioni comuni. Si potrebbe anche invece ad esse sostituirne altre. Io dubito che la spiegazione <u>filosofica</u> dell'origine sperimentale delle nozioni comuni, e dei tre postulati 10, 11 e 12 basti per renderle più affini agli assiomi che ai postulati.

Il Dini ti saluta e ti manda il ritratto del Novi. Rimanda sollecitamente le prove di stampa. Io sarò il dì 29 a Firenze. Potresti indirizzarmi le lettere là scrivendomi il 28. Poni sopra <u>Via San Gallo</u> N.° 33. 2°. P°. Dimmi anche se crederesti conveniente di stereotipare, e tirare gli esemplari pochi alla volta. Il 5 Ottobre dovrò essere nuovamente a Firenze per una nuova Commissione per esaminare le proposte della Facoltà Universitaria. Spero di trovarci Brioschi, al quale farò leggere quello che ho scritto sopra l'insegnamento dell'aritmetica e dell'Algebra nei Licei, dove mi vien fatto di toccare anche di volo i vantaggi dell'Euclide.[36]

---

[33] Cremona, Betti e Brioschi stanno lavorando al volume: E. Betti, F. Brioschi, *Gli elementi di Euclide con note aggiunte ed esercizi ad uso de' Ginnasi e de' licei*, Firenze, Le Monnier, 1867 e questo sarà l'argomento preminente nelle lettere seguenti.

[34] Si tratta probabilmente di H. Hankel, "Untersuchungen über die unendlich oft oszillirenden und unstetigen Funktionen" del 1870.

[35] J. Hoüel, *Essai critique sur les principes fondamentaux de la Géométrie élémentaire*, Paris, 1867.

[36] In [Menghini, 1996], lettere nn. 46 e 48, pp. 37 e 39, Cremona spingeva Betti a scrivere un "opuscolo sull'insegnamento dell'aritmetica e dell'algebra". Non si ha notizia che tale libro sia stato effettivamente pubblicato, ma, dal contenuto di questa lettera, è evidente che Betti stesse preparando del materiale e anche che si stesse confrontando sia con Brioschi (il quale era Membro del Consiglio Superiore della Pubblica Istruzione) che con Cremona

Ama il

Tuo Amico aff
Enrico Betti

**23**

Pistoia, 30 Settembre 1867

Mio caro Cremona,

Ho ordinato che si stereotipi e si tirino per ora 3000 Copie. per la seconda tiratura si potranno fare piccole correzioni sopra le lastre, e si aggiungeranno in fine quelle Note che crederai conveniente. Mi pare che Euclide debba rimanere il testo come in Inghilterra. Il Colenso[37] della Scuola Normale dove rimase? non lo prendesti tu? Napoli mi disse che i programmi sono sotto i torchi. Dovendo essere testo per le scuole credo che andranno subito le 3000 Copie e altre e vi sarà un sufficiente guadagno da dividere. Mandami con precisione il titolo, per annunziarlo. *[...]*: ad uso dei Licei, per cura dei Professori Brioschi e Betti.

Dini pochi giorni sono aveva ricevuto e rivisto le prove di stampa che ti avrà rimandato.

Procura di lavorare anche tu alla edizione di Euclide, perché non si debba *[...]* troppo. La fretta può essere una scusa, ma non bisogna trovarcisi troppo.

Scrivimi presto ed ama sempre

il tuo Amico aff
Enrico Betti

**24**

Mio caro Cremona

Se mi rimandi per qualche giorno Colenso[38] tradurrò io, o farò tradurre da qualcuno gli esercizi. Quanto alla proposizione F, trovo che la tua osservazione sulla dimostrazione è giusta. In Euclide non si trova fatto uso di ciò che non si è dimostrato che si può costruire, o che è postulato o domandato di poter costruire. Ora soltanto nel Libro VI, Euclide dimostra la possibilità di determinare una parte aliquota qualunque di una data linea terminata. Si potrebbe anche la proposizione F dimostrare in altro modo, fondandosi sopra la proposizione 17 e dimostrando che se la 1° e la 2° hanno una misura comune, l'hanno anche la 3° e la 4°, e quindi che la 3° e la 4° sono egualmente moltiplici di questa misura, come lo sono la 1° e la 2° della loro. Ma converrà aggiungere tutto questo? Forse negli altri Libri degli Elementi che non sono in Viviani[39] vi sarà qualche cosa. Poiché la Proposizione E senza la F mi pare che non sia completa, le lascerei tutte due. Intanto non le fo comporre come pure le proposizioni L ed M.

Ama sempre

il tuo amico aff
Enrico Betti

Pistoia, 24 Ottobre 67.

---

su tale argomento. Probabilmente il materiale prodotto venne utilizzato per la stesura del documento "Istruzioni e programmi per l'insegnamento secondario classico e tecnico, normale e magistrale, ed elementare nelle pubbliche scuole del Regno" pubblicato sul 1° Supplemento al N. 291 della *Gazzetta Ufficiale del Regno d'Italia* (24 ottobre 1867), al paragrafo "Istruzioni e programmi per l'insegnamento della matematica nei ginnasi e nei licei", dove, ad esempio, si legge "Nella geometria, per dare all'insegnamento la massima efficacia educativa, e per ridurre a un tempo la materia entro modesti confini, basta applicare alle nostre l'esempio delle scuole inglesi, facendo ricorso agli *elementi di Euclide*, che per consenso universale sono il più perfetto modello di rigore geometrico" e dove i libri di Euclide sono essenzialmente i programmi dei diversi anni scolastici.

[37] Si tratta probabilmente di J.W. Colenso, *Euclid's Geometrical problems, as given in the Edition of Euclid's Elements*, Longman, London, 1846.
[38] Si veda la nota precedente.
[39] V. Viviani, *Elementi piani, e solidi d'Euclide ecc.*, in Firenze, per il Carlieri, 1690.

**25**
*Carta intestata:* R. Scuola Normale Superiore

Pisa 21 Novembre 1867

Mio caro Cremona,

Due sere fa vidi uno della Società Lemonnier il quale era entusiasmato dello spaccio dell'Euclide. Pare che nonostante la poca premura di divulgare la pubblicazione le fatiche frutteranno bene.

L'appendice sulle misure mi pare che converrà metterla in fine. Se non puoi penserò io.

Quanto alle figure farò osservazione quando vo a Firenze. Intanto tu guarda quelle del 5° Libro nelle bozze che ti ho mandato. Io aveva consigliato di fare come nel Peyrard,[40] le divisioni con una piccola interruzione invece che con una lineetta. Non so perché dopo averne convenuto hanno messo le lineette. Che te ne pare. Scrivigliene francamente il tuo parere.

Quanto alla Prefazione che è venuta così in ritardo quando mi pareva premesse far presto a pubblicare i primi libri già pronti perché si riaprivano le scuole, io che credeva fosse stata rivista da te, la guardai appena. Ma si può correggere. Mandami le correzioni; faremo ritirare la lastra di quella pagina e nelle Copie avvenire sarà corretta.

I problemi *[e le]* proposizioni da dimostrarsi che mi hai mandato non sono nel Colenso.[41] Saranno aggiunte.

Cercherò io del Corso di Matematiche ad uso degli aspiranti etc, e se lo trovo come spero, dal Menozzi che è librajo a Pisa e a Modena, lo comprerò per la Società etc.[42]

Mi pare che converrà non limitarsi alla misura dei solidi interi, ma estendersi anche alle parti di essi.

Dini ti saluta. Sarà pubblicata una sua Memoria sopra a vari punti della teoria della superfici nelle Memorie della Società Italiana.[43] Mi pare assai importante. È stato nominato Professore straordinario ma di ambedue le Cattedre di Geodesia e di Algebra, e quel che è peggio con sole 2000 Lire di stipendio, mentre l'anno scorso ne ebbe 3000. È stata una indegnità dovuta al solito Gatti. Io non ho potuto impedirla perché il Napoli mi diceva che avrebbero fatto quel che proponeva il Rettore,[44] e il Rettore proponeva di farlo straordinario di Algebra e di dargli l'incarico di Geodesia, con che lo stipendio sarebbe stato non vergognoso.

Il mio lavoro sull'insegnamento dell'Aritmetica e dell'Algebra nei Licei è rimasto alle 12 pagine che aveva scritto in campagna nel Settembre. Ora ricevo grandi sollecitazioni per fare un trattato di Aritmetica e di Algebra nei Licei. Intanto sino alla fine di questa settimana le occupazioni della scuola e gli esami non mi lasciano luogo a pensare a questo.

Bertini ha presentato la sua Tesi, e dimani avrà la discussione. Mi pare un buon lavoro. In quest'anno seguiterà, come Alunno di questa Scuola i suoi studi di Analisi, e poi verrà con te a seguitare a occuparsi di Geometria.

I programmi di Fisica e le Istruzioni mi pare che potessero esser fatti meglio. Felici pubblicherà nell'anno il suo trattato pei Licei.[45]

Salutami tanto Brioschi e Casorati, scrivimi presto ed ama sempre

il tuo Amico aff
Enrico Betti

---

[40] F. Peyrard, *Eleméns de géométrie d'Euclide traduits littéralement, et suivis d'un Traité du Cercle, du Cylindre, du Cône et de la Sphère; de la mesure des Surfaces et des Solides; avec des Notes; par F. Peyrard Bibliothécaire de l'École Polytechnique*, Paris, 1804 (disponibile in rete).
[41] Si veda la nota 37.
[42] P. Cassiani, *Corso di Matematiche ad uso degli aspiranti alla scuola di Artiglieria e genio di Modena*, Modena, 1806.
[43] U. Dini, "Sopra alcuni punti della teoria delle superfici", *Mem. d. Soc. ital. d. scienze* (XL), s. 3, 1868, pp. 17-92.
[44] Fausto Mazzuoli, rettore dell'Università di Pisa dal 1865 al 1870 e dal 1879 al 1880.
[45] Felici pubblicò numerosi articoli di fisica, ma l'unico volume da lui pubblicato pare essere *Appunti per lezioni di Fisica sperimentale*, Pieraccini, Pisa, 1884, quindi di molto posteriore a questa lettera.

**26**

Mio caro Cremona,

Ti porterà questo mio biglietto il Sig$^r$ Gino Della Rocca[46] che viene costì a far il corso dell'Istituto Tecnico Superiore. Egli prese la Licenza in Pisa l'anno decorso. Studiò la Geodesia teoretica sotto il Dini, il quale me ne parlò sempre come di uno dei suoi più bravi Scolari. Io te lo raccomando certo come sono che Egli si distinguerà tra i migliori costì come quà.
Rimanda sollecitamente le prove, perché altrimenti non si può andare avanti, e alla fine dell'anno si vogliono dare in luce altri tre libri.
Scrivimi ed ama sempre

<div align="right">il tuo Amico<br>Enrico Betti</div>

Pisa, 11 Dicembre 1867.

**27**
*Carta intestata:* R. Scuola Normale Superiore

<div align="right">Pisa 27 Dicembre 1867</div>

Mio caro Cremona,

Hai disponibile un poco di posto del prossimo fascicolo degli Annali? Poche pagine basterebbero. Nel trattare la teoria della propagazione del calore ho effettuato una integrazione con metodo analogo a quello tenuto da Riemann nella sua memoria sulla propagazione delle onde piane, che pubblicherei volentieri.[47]
Tra oggi e dimani lo stampatore avrà stampato il 4° libro e una parte del 5° e impaginato tutto il resto del 5°, il 6° e gli esercizi. Ti manderò le ultime prove che tu rimanderai subito, e al principio dell'anno saranno pubblicati questi altri tre libri. Seguiteranno subito al libro 11°: tu manderai il 12°. Io ordinai il Corso di Modena,[48] ma non l'ho anche avuto. Il Law[49] l'ho avuto, e l'ho dato allo stampatore. Del 1° Libro furono esauriti 5000 esemplari e tirati altri duemila, del 2° e 3°, ne furono esaurite pure un gran numero di copie, ma assai minore. Quando sarà stampato tutto il volume faremo il conto del guadagno netto, e verremo alla divisione.
Ho fatto mandare una copia a Platner ed una ad Houel, il quale scrisse a Forti che approvava moltissimo la nostra pubblicazione. Però a quanto pare egli crede che si fosse dovuto tradurre Euclide in linguaggio moderno come ha fatto Egli del 1° Libro.
Io sono incerto quanto alle due proposizioni di Simson[50] che sono l'ultime del 5° Libro. Il Law le mette: ma mi sembrano assai complicate.
Salutami Casorati e Brioschi quando tornerà da Firenze dove a quanto mi disse deve essere.
Ama sempre

<div align="right">il tuo aff Amico<br>Enrico Betti</div>

---

[46] Si veda [Menghini, 1996], lettera n. 51, p. 41.
[47] E. Betti, "Sopra la determinazione delle temperature variabili di una lastra terminata", *Annali di matematica*, s. 2, v. I, f. IV, 1867, pp. 373-380.
[48] Si veda la lettera precedente.
[49] H. Law, *Euclid's Elements*, London, 1855. Nel 1853 Law aveva pubblicato in due volume anche *Euclid with propositions, and notes, and an essay on logic*.
[50] Si tratta di una delle numerose ristampe dell'opera di Robert Simson su Euclide, la cui prima pubblicazione fu a Glasgow nel 1756.

**28**
*Carta intestata:* R. Scuola Normale Superiore

Pisa 14 Febbrajo 1868

Mio caro Cremona,

Il ritardo nella stampa di Euclide Lemonnier lo attribuisce al costruttore delle figure, che è stato ammalato. Queste figure di solidi richiedono molto più tempo perché si fanno come quelle del Law. Le prime due che vidi, erano quasi uguali a quelle. Speriamo ora che è guarito questo costruttore si andrà speditamente. Ti farò mandare le bozze di tutto il libro e in doppio esemplare. Il Libro XII lo ricevetti e lo passai a Lemonnier.

Il Peyron[51] è in Magliabechiana a Firenze e credo anche qui in Biblioteca: se non vi fosse lo procurerò, e lo vedrò. Mi par molto probabile che ce ne potremo valere.

Prepara pure subito tutte le correzioni da farsi al Libro I, e mandale appena le hai in pronti. Io ho già avvertiti gli Editori che in seguito avremmo fatto correzioni da richiedere anche di rifare intere pagine. Non so se siano ancora esaurite le ultime 2000 Copie che tirarono del Libro I, e sarebbe bene quando ne tireranno altre 2000 le tirassero colle correzioni, e che quando si pubblicheranno esemplari di tutta l'Opera, vi fosse nella medesima il primo libro corretto. Di ciò ne parlerò quando andrò a Firenze, ma ti posso assicurare sin d'ora che sarà fatto ciò che proponi, perché ne sono già prevenuti.

Salutami Brioschi, scrivimi ed ama

il Tuo Amico
Betti

**29**
*Carta intestata:* R. Scuola Normale Superiore

Pisa 7 Marzo 1868

Mio caro Cremona,

È composto tutto l'undecimo libro di Euclide, e Lunedì si comincia il $12^{mo}$. Martedì prossimo tu riceverai tutte le prove di stampa del L. $10^{mo}$ coll'originale. Parlai con Lemonnier delle correzioni al 1° Libro, e non vi è difficoltà. Però bisognerebbe che tu le mandassi presto, perché non si debbano tirarne altre copie prima che siano fatte.

Nel Colenso[52] vi sono pochissimi esercizi di Geometria solida, messi insieme con esercizi di Geometria spettanti alla teoria delle proporzioni. Che cosa dobbiamo fare? Bisognerebbe prenderne altrove. Quanto all'Appendice ora urge decidere. Il Peyrard di cui mi parlaste è quella edizione in tre volumi che è alla Magliabechiana?[53]

Tra pochi giorni avrai un Articolo del Dini per gli Annali.[54] Il mio quando sarà pubblicato?[55] Io devo tornare a Firenze il 14, e allora bisognerà fissare tutto colla stamperia purché l'Euclide intero possa essere pubblicato nella prima metà di Aprile.

Scrivimi ed ama

il Tuo aff Amico
Enrico Betti

---

[51] Probabilmente si tratta del Peyrard di cui alla nota 40.
[52] Si veda la nota 37.
[53] Il Peyrard è in un unico volume. Si veda [Menghini, 1996], lettera n. 54, p. 43.
[54] U. Dini, "Sulle superficie che hanno le linee di curvatura piane", *Annali di matematica*, s. 2, v. I, 1868, pp. 146-154.
[55] Si veda la lettera n. 19.

**30**
*Carta intestata:* R. Scuola Normale Superiore

Pisa 13 Aprile 1868

Mio caro Cremona,

Io ho pensato di pubblicare intanto i libri 11 e 12. L'appendice la pubblicheremo a sé e con tutto il libro. Però bisognerebbe sbrigarsi anche per quella. Si può prendere da Peyrard con quelle modificazioni ed aggiunte che crederemo. Intanto potresti farla subito tradurre. Pensa anche alle correzioni perché si possa mandar fuori tutto il libro in buono stato. È fatto il conto fino al 31 Dicembre, e delle vendite fatte fino a quel giorno ci tocca di parte L. 394,96. Delle vendite posteriori a quel giorno non è anche fatto il conto.

Mi ha detto Brioschi che vai a fare una ispezione nelle antiche Provincie.[56] Dove debbo indirizzarti le lettere? Quando avrò le prove di stampa dell'articolo per gli Annali?

Scrivimi ed ama

il tuo aff Amico
Enrico Betti

**31**
*Carta intestata:* R. Scuola Normale Superiore

Pisa 2 Maggio 1868

Mio caro Cremona,

Ti rimando subito le prove corrette. Fammi il piacere di rivedere le correzioni. Nella 2° pagina avevamo lasciato la divisione dal primo al 2° paragrafo. Gl'integrali definiti hanno i limiti troppo discosti. Poi è stato scritto nelle formule lim per limite collo stesso carattere delle lettere che denotano quantità. Non sarebbe bene scrivere lim. con altro carattere?[57]

Io feci subito ricerca del libro di cui mi scrivesti. Qui a Pisa quel Magazzino termina al 1806.[58] Ne scrissi a un impiegato della Magliabechiana a Firenze e quegli mi mandò l'anno 1805 dicendo che non aveva trovato altro numero in cui si parlasse dell'Ocia Perelliana. Ma in quello non vi erano che poche parole. Stasera vo a Firenze e ne farò ricerca da me.

Ama sempre

il tuo Amico aff
Enrico Betti

---

[56] Dal giugno Cremona stava facendo un giro di ispezione delle scuole del Piemonte e della Liguria su incarico del Ministero. Già nel 1863 aveva avuto questo incarico per le Puglie e le Calabrie. Si vedano [Menghini, 1996], lettere n. 8 e seguenti, da p. 13 e la n. 57, p. 47. Si veda anche qui la lettera n. 5.
[57] Betti sta correggendo le bozze del proprio articolo: si veda la lettera n. 27.
[58] In [Menghini, 1996], lettera n. 55, p. 45, Cremona chiedeva a Betti di cercare un opuscolo di Pietro Ferroni, dal titolo *Ocia perelliana*, che, secondo le sue informazioni, avrebbe dovuto essere contenuto nel *Magazzino di letteratura, scienze, arti, economia politica e commercio*. In realtà qui si trova soltanto un "Ragguaglio" del manoscritto di Ferroni a cura di P. Pozzetti (v. 1, 1805, pp. 150-156; consultabile in rete).

**32**
*Carta intestata:* R. Scuola Normale Superiore

Pisa 18 Luglio 1868

Mio caro Cremona,

Quanto alla Presidenza della Società Italiana ne aveva già scritto a Genocchi e parlato con questi di Pisa. Genocchi, con i Sismonda e Moris è d'accordo con noi, e quelli di Pisa ugualmente. Credo che la elezione riuscirà a seconda dei nostri desideri.[59]

Quanto al Baltzer io ho già fatto il rapporto.[60] Non ho proposto l'approvazione, perché dietro le norme stabilite per questa non si poteva perché non corrisponde pienamente, essendo fatto per delle scuole dove la Matematica s'insegna più estesamente che nei nostri Licei. Ma del resto ne ho riconosciuto il valore, e ho detto che nessuno insegnante può mancare di un tal libro e nessuno dei giovani che intende seguitare gli studi matematici alle Università, ed ho detto che tu e gli Editori avete reso un vero servigio all'insegnamento della Matematica con quella traduzione.[61]

Il 1° Libro è stato incorporato e ristampato. Io credeva che le ultime correzioni e aggiunte le avessi fatte tu, e non ho fatto che completarle con i 4 trattati inglesi che ho, e le ho trovate conformi in generale, e conformi anche ad alcune osservazioni che mi comunicò Brioschi del Sacchi, e che aveva avute da altri Professori di Liceo.

Quanto a questa nostra pubblicazione bisogna pensare alla divisione del guadagno: bisogna stabilire la proporzione in cui dev'essere diviso. Giustizia vuole che sia diviso fra te, me e Brioschi; tu hai lavorato più di tutti; io ho lavorato e messo il nome, Brioschi ha messo il nome e fatto la prefazione. Si potrebbe fare 2, 2, 1? Se questa proporzione non ti par conveniente dimmelo francamente ed io farò come tu credi meglio.

Io ho veduto Brioschi a Firenze ieri e lo rivedrò il 26 di questo mese e parlerò di quello che mi dici relativamente alla Giunta. Da Casorati ebbi lettere 8 o 10 giorni fa.

Ti saluta il Dini. Teza partì di qua Martedì prossimo passato.
Quanto ai corsi di Aritmetica e Algebra aveva incaricato il Finzi di occuparsene insieme con me. Ma per ora si è fatto poco. Penserei intanto di rifare alcuni Capitoli del Bertrand,[62] che cosa ne dici?

Scrivimi presto ed ama il

Tuo Amico
Enrico Betti

---

[59] Betti e Cremona desideravano portare Brioschi alla Presidenza, cosa che avvenne. Si vedano [Menghini, 1996], lettere nn. 57 e 58, pp. 47-49.

[60] Cremona stava curando la traduzione del libro R. Baltzer, *Die Elemente der Mathematik*, Leipzig, 1865 che uscì con il titolo *Elementi di Matematica. Prima versione italiana fatta sulla seconda edizione di Lipsia*, Genova, 1865-1867. Nella lettera n. 56 ([Menghini, 1996], p. 46) del 2 maggio, Cremona aveva chiesto espressamente che "il Consiglio Superiore dichiari l'opera utile per le scuole", cosa che Betti non fece. Cremona si dimostrò molto contrariato da questa posizione che gli faceva anche perdere una certa entrata finanziaria. Si vedano [Menghini, 1996], lettere nn. 56-61, pp. 46-51. Si vedano anche qui le lettere seguenti; in particolare la 35 nella quale Betti trascrive la sua relazione al Consiglio.

[61] Anche nel libro di Euclide (si veda la nota 33) Betti e Brioschi scrivono a proposito della traduzione del Baltzer: "Raccomandiamo quest'opera a quei giovani che intendono progredire negli studi matematici".

[62] J. Bertrand, *Traité élémentaire d'algébre*, Paris, Hachette, 1850. Betti nel 1856 aveva pubblicato il *Trattato di Algebra elementare di Giuseppe Bertrand, Prima trad. it. con Note ed Aggiunte* (Firenze, LeMonnier) proponendone poi "parecchie ristampe successive", come si legge nell'elenco dei lavori scientifici a p. VIII delle *Opere Matematiche di Enrico Betti*.

**33**
*Carta intestata:* R. Scuola Normale Superiore

Pisa 23 Luglio 1868

Caro Cremona,

Se si fosse trattato di una somma per una sola volta io avrei certamente proposto che si desse interamente a che ti aveva aiutato;[63] ma trattandosi invece di una vendita annua non piccola il 1° anno, più grande nei successivi, perché non occorre più la spesa della composizione della stereotipia, e che durerà finché Euclide sarà testo nelle nostre scuole, questo non mi parve e non mi pare ragionevole.

Cominciamo dunque da prelevare il compenso da darsi a Platner che tu stabilirai, quanto al resto se non vuoi nulla tu, intenderai facilmente che non potremo a più forte ragione prender nulla io e Brioschi. Quindi o si divide in parti uguali o in altro modo che paia più conveniente o se ne dispone ad un altro uso tutti d'accordo. Quanto a me certamente non si troveranno difficoltà. Intanto subito vo a Firenze e là chiederò il conto del nostro credito; ve lo farò sapere, e stabiliremo il da farsi tra noi, come buoni amici, pieni di fiducia reciproca, e che sopra ogni altra cosa hanno a cuore l'utile degli studi e la propria dignità.

Quanto al Baltzer tu non hai conoscenza, pare, delle norme stabilite dal Consiglio. Egli non approva altro che i libri che possono con vantaggio straordinario essere adottati come libri di testo per le nostre scuole. Ora se tu avessi dovuto fare il rapporto io sono convintissimo che lo avresti fatto precisamente come l'ho fatto io. Del resto se mi sono ingannato ed è possibile fare altrimenti, avresti dovuto essere certo che niente mi sarà più gradito.[64]

I saluti da parte del Dini.

il Tuo Enrico Betti

**34**
*Carta intestata:* R. Scuola Normale Superiore

Pisa 22 Agosto 1868

Amico carissimo,

Mi scrive Lemonnier che parecchi corrispondenti gli hanno ripetutamente chiesto il motivo dell'indugio alla pubblicazione dell'Appendice sulle misure promessa nella Prefazione.

Tu mi scrivesti molto tempo indietro di avere incaricato Platner di farne la traduzione dal Peyrard. Mi manifestavi anche il desiderio che avresti avuto di sopprimerla se fosse stato possibile. Ma contro questo desiderio vi sono due ostacoli: la promessa fatta nella Prefazione e quel che più importa la esigenza dei programmi.

Scrivimi qualche cosa in proposito e dimmi ancora dove debbo fare indirizzare le Lire 250 per Platner. Quanto all'altra somma aspetto che decidiamo qualche cosa relativamente al premio. Ne parleremo con Brioschi nel Settembre.

Tu hai occasione di venire presto a Firenze? Nel caso affermativo, fammelo sapere; io il primo di Settembre vo in campagna a Pistoia e verrei a trovarti a Firenze per istare un poco insieme e parlare di molte cose.

Neumann e Clebsch pubblicheranno questo nuovo Giornale,[65] che stamperà le memorie anche in lingua Italiana. Mi piace molto questo affratellamento nella Scienza.

Scrivimi presto ed ama sempre

il Tuo Amico
Enrico Betti

---

[63] Giacomo Platner aveva curato le traduzioni per l'Euclide e Cremona aveva chiesto (si veda [Menghini, 1996], lettera n. 58, pp. 48) che la somma destinata a lui fosse versata direttamente all'amico.
[64] Si veda la lettera precedente.
[65] *Mathematische Annalen, herausgegeben von A. Clebsch und C. Neumann*, Berlin. Giornale fondato nel 1868 e chiuso tra il 1945 e il 1947.

**35**

Pistoia, 18 Settembre 1868

Mio caro Cremona,

Feci spedire L. 250 a Platner. L'ultima volta che parlai a Lemonnier non avevano avuto avviso che le avesse ricevute. Importerebbe però moltissimo per liquidare i nostri conti, consegnare la traduzione per l'Appendice, e quindi se Platner non intende farla, scrivimene subito perché penserò io a farla eseguire da altri. Se non ha difficoltà a farla, sollecitalo a sbrigarsi.

Il mio rapporto sopra il Baltzer come ti avrà detto Brioschi non sarà pubblicato, soltanto sarà comunicato all'editore, il quale potrà farne quell'uso che crederà, potrà anche stamparne uno per te o scritto in una prefazione. Io te lo trascrivo qui come mi hai chiesto.[66]

Gli Elementi di Matematica di Baltzer tradotti dal Cremona sono libri eccellenti ed hanno fatto opera utile all'insegnamento in Italia il traduttore e l'editore col renderlo noto tra noi; ma questi libri sono fatti per iscuole dove si dà alla Matematica elementare una attenzione di gran lunga maggiore a quella che le si dà nelle nostre Scuole Secondarie. Però mentre il vostro Relatore crede che ogni insegnante di Matematica dei nostri Ginnasi Licei e Scuole Tecniche dovrebbe prenderne cognizione, e che quei giovani i quali hanno finito il Corso di Matematica nei Licei e si accingono ad intraprendere gli Studi di Matematica nelle Università studiandoli ne ricaveranno un ottima preparazione agli Studi Universitari, non vi può proporre che siano approvati come libri di testo per le nostre Scuole Secondarie.

La Conclusione del rapporto e il rapporto stesso furono approvati dalla Commissione, però questo non impedirà qualche mutamento nel rapporto stesso, e qualche aggiunta o soppressione che lasciando ferma la conclusione possa essere più utile all'Editore. Se tu me la indichi io alla nuova riunione del Consiglio farò restituire la pratica e ve la introdurrò.

Di libri approvati non ve ne sarà alcuno; bisognerà pensare a prendere anche qualche misura relativamente ai libri approvati dal Consiglio precedente e dai Comitati, perché la misura è molto differente. L'approvazione quale prima si richiedeva dal Consiglio Superiore è data ora dai Consigli provinciali scolastici. Quella del Consiglio Superiore, insieme colla radiazione *[...]* dalle liste dei Consigli provinciali, è riservata per i soli trattati che corrispondono in modo straordinario a tutti i bisogni dell'insegnamento cui si riferiscono.

Scrivimi ed ama sempre

il tuo aff Amico
Enrico Betti

---

[66] Si veda la lettera n. 32.

**36**

Pistoia, 2 Ottobre 1868

Mio caro Cremona,

Non posso trascriverti il rapporto sopra il libretto del Zinna,[67] perché fu passato alla Divisione. Però posso dirti che non era favorevole. Quel che vi era di nuovo era erroneo; e come Monografia non esponeva lo stato attuale della Scienza, e questo in sostanza era quello che io diceva.

Il Parlatore s'incaricò di scrivere e di mandare le schede,[68] ed io credo aver fatto molto male ad affidarmi a Lui. Dici a Stoppani che mandi subito a me lo scritto sopra un foglio di carta chi egli elegge a presidente e lo firmi. È passato il tempo ed io avrei voluto far subito lo spoglio delle schede, ma il Sig. Parlatore non si sa dove sia, e a me tocca ad aspettare e ad essere negligente come Lui. La tua la ricevetti in tempo debito.

Ho già mandato a Lemonnier il manoscritto. Ho mutato un poco in principio. Vi era definita la misura per mezzo del rapporto, e poi spiegato il significato del rapporto per mezzo della misura. Ecco come io ho fatto. Non riferisco le parole precise perché non ho qui il manoscritto.

Date due grandezze omogenee A e B se la moltiplica di B secondo un numero $\underline{m}$ è uguale ad A, il numero $\underline{m}$ si dice la misura di A, quando si prende per unità B.

Se non vi è moltiplica di B uguale ad A, ma la moltiplica di B secondo il numero $\underline{m}$ è uguale alla moltiplica di A secondo il numero $\underline{n}$, cioè se la moltiplica della $\underline{n}^{esima}$ parte di B secondo il numero $\underline{m}$ è uguale ad A, il numero $\frac{m}{n}$ si dice la misura di A quando si prende per unità B. Se poi non vi sono moltiplici di A uguali a moltiplici di B, prendendo la moltiplica di A secondo un qualunque numero $\underline{n}$, si potrà trovare $\underline{m}$ tale che la moltiplica di B secondo m sia minore e la moltiplica secondo m+1 sia maggiore di nA e quindi:

$$\frac{m}{m}B < A < \frac{m+1}{n}B$$

Quindi $\frac{m}{n}$ ed $\frac{m+1}{n}$ sono le misure di due grandezze una maggiore e una minore di A, e differenti da A meno della $n^{esima}$ parte di B; il limite di questi numeri che è unico e determinato dà la grandezza di A quando è noto B; quindi può definirsi la misura di A quando è preso per unità B. Questi limiti non sono numeri interi né frazionari, ma in quanto esprimono le misure di grandezze, ci si possono fare sopra i medesimi, come si dimostra in aritmetica, operazioni che seguono le stesse leggi di quelle che si fanno sopra i numeri interi e frazionari, si chiamano numeri anch'essi e si dicono incommensurabili.

In tutti e tre i casi il numero che esprime la misura di una grandezza A quando è presa per unità un altra grandezza B, è uguale (v. Baltzer) al rapporto $\frac{A}{B}$ dei numeri che esprimono le misure di A e B quando si prenda per unità una grandezza arbitraria; perciò si dice rapporto delle grandezze A e B.

Quindi viene il teorema che stabilisce che se due *[regioni]* nel senso di Euclide sono uguali, saranno uguali anche i numeri che ne esprimono i rapporti.

Questo punto è assai delicato e vorrei che fosse chiaro ed esatto. Però ho voluto scrivertelo perché tu mi faccia le tue osservazioni. Ama sempre

il tuo Amico
Enrico Betti

---

[67] A. Zinna, *Monografia sulla risoluzione delle equazioni numeriche*. Si veda [Menghini, 1996], lettera n. 58, pp. 48.
[68] Si tratta delle schede per la nomina del Presidente della Società Italiana di Scienze (detta dei XL). Cremona, nella lettera del 30 settembre (si veda [Menghini, 1996]) aveva scritto a Betti di aver votato per Brioschi che in effetti sarà eletto e rimarrà in carica fino al 1874.

**37**
*Carta intestata:* R. Scuola Normale Superiore

Pisa 17 Novembre 1869

Mio caro Cremona,

Ti mando la tesi del <u>Levi</u>. Questa dissertazione la conosceva prima che la scrivesse, perché non è altro che una copia con qualche mutamento di parole di una parte delle lezioni che io detti a Pisa (quando vi era scolaro il Levi) sopra la teoria della elasticità. Dunque <u>non vi è merito alcuno di originalità, e non prova anche a mostrare la sua cultura</u>. Però mi rammento che Egli era un giovane molto intelligente, e studiò assai seguendo i miei corsi, ma anche da scolaro era inferiore a Bertini.[69] È di natura piccolissimo e storpio.

Salutami Brioschi e Casorati, ricevi i saluti del <u>Dini</u>, che ha già pronto il suo lavoro.[70] Ama sempre
il tuo amico
Enrico Betti

**38**
*Carta intestata:* R. Scuola Normale Superiore

Pisa 30 Aprile 1870

Mio caro Cremona,

Eccoti la copia dell'autografo che Jacobi aveva dato a Lavagna nel 1843 al Congresso di Lucca nel quale Lavagna era Segretario della Sezione delle Scienze Fisiche e Matematiche. Questa proposizione si trova nella Memoria sopra il Moltiplicatore §12.

Quanto alla Memoria del Volpicelli,[71] poiché il teorema, già molto vecchio, è vero, non credo ne metta il conto a perdere il tempo a criticare le dimostrazioni del Volpicelli.

Ti avrà detto Brioschi che degli utili del 69 pagheranno L. 250 il mese a cominciare dagli ultimi di Marzo. Fissammo di riscuotere un mese per ciascuno. Il mese di Marzo lo riscossi io, quello di Aprile sarà pagato a te, quello di Maggio a Brioschi e così di seguito fino alla fine dell'anno. La somma dell'ultimo mese sarà divisa in parti uguali.

Salutami Casorati ed ama
il tuo aff Amico
Enrico Betti

---

[69] Simeone Levi. Si veda [Menghini, 1996], lettera n. 65, p. 53 nella quale Cremona chiede a Betti un giudizio sulla dissertazione dal titolo "Sull'equilibrio di un corpo elastico" di Levi. Cremona era nella Commissione per il concorso alla cattedra di matematica nel liceo Parini, per la quale concorreva anche Bertini e sperava che questi non venisse superato.
[70] Potrebbe trattarsi di: U. Dini, "Sui prodotti infiniti", *Annali di matematica*, s. 2, v. II, 1869, pp. 28-38.
[71] P. Volpicelli, "Della distribuzione elettrica dei conduttori isolati", *Annali di matematica*, s. 2, v. III, 1870, pp. 249-268. Si veda [Menghini, 1996], lettera n. 66, p. 54.

**39**

Ponsacco presso Pontedera
18 Agosto 1870

Mio caro Cremona,

Ieri seppi dal Cassiere della Società Lemonnier che ti aveva spedito le 250 Lire soltanto il 12 di questo mese. Io gli aveva lasciato un biglietto sul suo tavolino fino dal 28 del mese di Luglio col quale lo avvisava che ti facesse questa spedizione: ma egli si è scusato col dirmi che il biglietto per vari giorni era rimasto tra i suoi fogli senza che lo avvertisse.

Io sto scrivendo un Articolo sopra la connessione degli spazi di un numero qualunque di dimensioni, e sopra la estensione agl'integrali di un numero qualunque di variabili del teorema relativo alla integrazione lungo una linea situata sopra una superficie più volte connessa, dal quale dipende la determinazione del numero dei moduli di periodicità delle funzioni che così si ottengono.

Pubblicherei volentieri questo articolo negli Annali, ma dubito che ci sia da aspettare molto tempo.[72] Fammi il piacere di scrivermi quando fai conto che possa cominciare la stampa del fascicolo in cui potresti inserire questo mio articolo. La lettera puoi indirizzarla <u>Pontedera</u> per <u>Ponsacco</u>.

Salutami Brioschi e Casorati ed ama sempre

il tuo Amico
Enrico Betti

**40**
*Carta intestata:* R. Scuola Normale Superiore

Pisa 3 Novembre 1870

Amico carissimo,

Ho ricevuto le prove di stampa, ma soltanto una parte; manca l'ultimo paragrafo nel quale applico il teorema sopra la connessione degli spazi agli integrali multipli e dimostro come serva alla determinazione della polidromia delle funzioni di più variabili alle quali danno origine. Ti prego a farmi il favore d'informarti come va questa cosa, e se non hanno composto quel paragrafo a farlo comporre perché sarebbe bene che fosse pubblicato l'articolo tutto completo in un solo fascicolo.

Detti ordine al Lemonnier che ti mandasse il primo di questo mese la tua parte utili Euclide *[sic!]*; il 6 andrò a Firenze e gli rammenterò la cosa nel caso che non te l'abbia mandata. Con questa sarà, colla differenza di poche lire, esaurito l'utile del 1869. Nell'anno venturo sapremo quanto sarà l'utile del 1870. Ricevi i saluti del Dini.

Saluta Casorati ed ama sempre

il tuo Amico
Enrico Betti

---

[72] E. Betti, "Sopra gli spazi di un numero qualunque di dimensioni", *Annali di matematica*, s. 2, v. IV, 1871, pp. 140-158.

**41**
*Carta intestata:* R. Scuola Normale Superiore

Pisa 25 Dicembre 1870

Caro Cremona,

Io non ho veduto la Memoria di Lipschitz;[73] mandamela ed io la tradurrò in italiano dentro un mese o un mese e mezzo.

Quanto ai programmi intorno ai quali mi facevi una dimanda nell'altra tua lettera, io ebbi una lettera del Ministero per proporre i cambiamenti che credessi opportuni. Io non risposi altro che a voce a Barberis,[74] e gli dissi che non avrei fatto altri mutamenti, fuori che avrei messo qualche cosa per il terzo anno del Liceo. Ho saputo poi che hanno fatto dei mutamenti e non tutti buoni a mio credere.

Ho ricevuto il saldo dell'Anno 1869 degli utili di Euclide. Ti mando i 10 franchi che ti pervengono di questa ultima riscossione. Non è ancora fatto il conto degli utili del 1870.

Salutami Casorati; tanti saluti dal Dini. Ricevi mille buoni auguri per il nuovo anno ed ama sempre

il tuo Amico aff
Enrico Betti

**42**
*Carta intestata:* R. Scuola Normale Superiore

Pisa 26 Marzo 1871

Caro Cremona,

Ho ricevuto le bozze di stampa della traduzione dell'Articolo di Lipschitz;[75] fammi il favore di rimandarmi per la posta l'originale tedesco. Dimani sera vo a Firenze, puoi indirizzarlo (se lo mandi subito) là, in Via San Gallo N.° 33. 2° Piano. Ho ricevuto il resoconto della Società Lemonnier. Siamo accreditati per L. 1795,57. Ce li pagheranno come l'anno decorso un tanto il mese. Fisserò in questi giorni coll'Amministratore e faremo come si fece nell'anno passato, cioè prenderemo un mese per ciascuno. Sono mi pare più di 200 Lire meno dell'anno precedente. Verrà a trovarti il D'Arcais a cui ho dato un mio biglietto. È un bravo giovine, ed è un ricco signore, quindi può coltivare liberamente la Scienza. Presentalo a Casorati; avrà piacere di parlare di funzioni Abeliane delle quali si è occupato.

Salutami Casorati, scrivimi ed ama sempre

il Tuo aff Amico
Enrico Betti

*[seguono appunti scritti da Cremona]*

---

[73] L'articolo originale è: R. Lipschitz, "Beiträge zur der Theorie der Umkehrung eines Functionensystems", *Nachrichten von der Königl. Gesellschaft der Wissenschaften und der Georg-Augusts-Universität zu Göttingen*, 1870, pp. 439-477; la traduzione è: "Sopra la teoria della inversione di un sistema di funzioni", *Annali di matematica*, s. 2, v. IV, 1871, pp. 239-259, non firmata.
[74] Potrebbe trattarsi del commendator Giuseppe Barberis.
[75] Si veda la lettera precedente.

**43**
*Carta intestata:* R. Scuola Normale Superiore Pisa

Pisa, li 16 Giugno 1871

Caro Amico,

Aveva stabilito di fare come l'anno passato, cioè di riscuoter la somma in 9 rate mensili a cominciare dall'ultimo di Aprile. La prima rata il 1° di maggio le presi io; per la seconda aveva dato l'ordine di spedirla a te il 1° di Giugno. Pare che se ne fossero dimenticati, ma appena ricevuta la lettera ho scritto al Cassiere che facesse subito la spedizione, e credo che l'avrà fatta. Il 1° di Luglio deve riscuotersi un altra rata, che sarebbe per Brioschi. Se Egli non ha difficoltà potrei farla rimettere a te, a Lui potrei dare quella di Agosto ed io mi potrei prendere invece quella di Settembre. Farò come meglio accomoderà.

Ti manderemo i resoconti per l'<u>Jahrbuch</u>.[76] Salutami Brioschi e Casorati. Il Dini ti saluta. Ama sempre

il tuo Amico aff
Enrico Betti

**44**
*[in alto, scritto da Cremona:* (Far mandare il ms. al Fontebasso)*]*[77]

Sojana, 23 Agosto 1871

Caro Cremona,

L'articolo del Professor Marsano non mi pare che sia da inserirsi né tutto né in parte negli Annali. Il soggetto sarebbe adattato per il Giornale del Battaglini; ma io non ho avuto la pazienza di leggerlo fino in fondo ed ho fatto più presto a ritrovare da me la dimostrazione. Quindi io declinerei l'offerta.

Dietro la insistenza del Ministero è bisognato che il Consiglio superiore faccia un programma di concorso per un libro di Matematica elementare da servire per le scuole secondarie classiche. Si è proposto un concorso per un trattato di Aritmetica, Algebra ed elementi di trigonometria seguendo le norme e le istruzioni dei programmi del 1867, e per un trattato di Geometria che nella Geometria piana sia l'Euclide rivisto e nella solida aggiunga sempre al metodo Euclideo anche i concetti introduttivi dopo Euclide. Io ho già riferito in questo tempo al Consiglio che ha approvato; ma ora devo fare la relazione da rimettersi al Ministro.[78] Se mi potessi mandare quell'articolo che scrivesti in risposta al Rubini e dirmi la tua opinione in proposito, mi faresti cosa gratissima.[79] Non so che cosa armeggino al Ministero, ma io dubito che seguitino per la medesima via.

Io ho ripensato ancora a quello di cui parlammo a Pisa relativamente ai sistemi 4 volte infiniti di rette cioè all'elemento lineare di questi spazi a quattro dimensioni. Prendendo il momento di due rette infinitamente vicine un corpo solido in questo spazio sarebbe l'insieme di tutte le rette le cui coordinate sodisfarebbero *[sic!]* a certe diseguaglianze. Questo corpo sarebbe rigido quando fossero invariabili i momenti di tutte le rette infinitamente vicine che ne fanno parte. Poiché la curvatura di un tale spazio sarebbe nulla, potrebbe sovrapporsi a sé stesso in ciascuna sua parte, e quindi un corpo rigido potrebbe muoversi. Il moto di un tal corpo rigido dipenderebbe da 10 sole quantità, mentre nello spazio ordinario a tre dimensioni dipende da 6. Fissando in esso tre rette rimane dunque possibile il moto che dipende allora da una sola quantità, e in questo moto rimangono fisse tutte le rette che rimangono in quelle tre sopra una stessa superficie di linee, e le altre linee rette si muovono lungo superficie di linee. Fissando 4 rette si fissa il corpo, a meno che questa rette non siano tutte sopra la stessa superficie di linee (Iperboloide). Le analogie sono curiose e complete e più se ne potranno trovare; ma dove si trova sempre una differenza è nella espressione analitica dell'elemento lineare. L'elemento lineare negli spazi a tre e a due dimensioni sono radici quadrate di forme differenziali di 2° grado omogenee sempre positive: mentre i momenti delle rette sono anche negativi e la forma differenziale è omogenea di 2° grado ma non è una forma positiva, quindi per prenderne la radice per fondamento delle misure si ha l'immaginario. Io andava pensando se non sarebbe questo uno dei casi accennati da Riemann in cui si debba prendere per elementi delle misure la radice quarta del quadrato del momento. Insomma il momento si debba nella misura come le distanze nello spazio prendere sempre positivo.

---

[76] In [Menghini, 1996], lettera n. 71, p. 56, Cremona chiedeva a Betti di scrivere qualche breve estratto dei propri scritti da mandare al *Jahrbuch über die Fortschritte der Mathematik* di Berlino in modo da far conoscere all'estero i lavori italiani.

[77] Si veda [Menghini, 1996], lettera n. 63, p. 52.

[78] Cesare Correnti.

[79] Betti si riferisce alla polemica intercorsa sulle pagine del *Giornale di Matematiche* a proposito dell'Euclide. Si veda [Menghini, 1996], lettera n. 64, p. 52. L'articolo è: F. Brioschi, L. Cremona, "Lettera al direttore del Giornale di Matematiche", *Giornale di Matematiche*, v. VII, 1869, pp. 51-54, pubblicato anche nel terzo volume delle *Opere* di Cremona [Cremona, *Opere*], n. 82, p. 133.

Quando mi scrivi puoi sempre indirizzarmi la lettera a Pisa, oppure a Peccioli per Sojana ma allora bisognerà aggiungervi anche Toscana o Provincia di Pisa.

Ho lasciato l'ordine al Cassiere di Lemonnier che ti spedisca le Lire 200 a Milano i primi di Settembre. Parlai coll'Amministratore del tuo Trattato.

Scrivimi ad ama

il tuo aff Amico
Enrico Betti

**45**
*Carta intestata:* R. Scuola Normale Superiore Pisa

Pisa, li 12 Novembre 1871

Caro Cremona,

Avrai ricevuto da Battaglini i documenti relativi al concorso per la Meccanica razionale a Pisa. Credo che li avrai esaminati, e ti prego a mandarli più sollecitamente che puoi a Beltrami, al quale dirai che li rimetta a me appena che li avrà esaminati. Allora stabilirò il giorno in cui ci dovremo trovare in Pisa per prendere le deliberazioni intorno alle proposte da farsi.

Salutami Brioschi e Casorati ed ama sempre

il tuo aff Amico
Enrico Betti

**46**
*Carta intestata:* R. Scuola Normale Superiore Pisa

Pisa, li 4 Dicembre 1871

Caro Cremona,

Il 7 è seduta di Consiglio superiore a Roma; quindi la nostra adunanza è differita. Il Dini ti saluta e ti prega a mandargli l'Hankel[80] sopra le teorie infinitamente spesso oscillanti, che te lo renderà quando vieni a Pisa; l'ha ordinato ma non l'ha potuto ancora ricevere. Il Padova vorrebbe tradurre la Geometria Descrittiva di Fiedler,[81] che cosa ne dici? Scrissi a Lemonnier che ti mandasse i 200 franchi.

Ama sempre il

Tuo Amico aff
Enrico Betti

*[sul retro disegni e appunti scritti da Cremona]*

---

[80] Si veda la lettera n. 22 e [Menghini, 1996], lettera n. 74, p. 58.
[81] *Trattato di geometria descrittiva del dr. Guglielmo Fiedler, tradotto dall'Ingeg. Antonio Sayno e dal dott. Ernesto Padova. Versione migliorata coi consigli e le osservazioni dell'Autore e liberamente eseguita per meglio adattarla all'insegnamento negli Istituti Tecnici del Regno d'Italia*, Firenze, Le Monnier, 1874. Si veda la nota precedente.

**47**
*Carta intestata:* R. Scuola Normale Superiore Pisa

Pisa, li 5 Maggio 1872

Mio caro Cremona,

Quanto all'epoca della riunione della Commissione per il concorso a Torino non fisserò niente senza andarne d'accordo con Te. Mi pare che si potrebbe fare nella seconda metà di Giugno. Prima forse non sarà stato finito l'esame dei titoli di tutti; dopo vi sono gli esami, poi le bagnature e le villeggiature. Intanto appena avrò ricevuto alcuni documenti che mancano comincerò dal mandare a Te e al Beltrami l'incartamento, perché possiate esaminarle i titoli prima di partire per le vostre ispezioni.

La causa del ritardo a rimborsarci delle spese incorse per il concorso di Pisa sai tu qual è stata? l'aver terminato il fondo stabilito in Bilancio per questo oggetto, prima della fine dell'anno: credo che vi avranno rimediato in qualche modo. Io debbo tornare a Roma il 25 di questo mese, e allora farò nuove premure perché sbrighino la cosa.

Quest'anno gli utili Euclide sono anche minori; sono circa mille franchi. Dimani vo a Firenze e fisserò il modo del pagamento che sarà credo come nell'anno passato.

Nel Giugno io dovrò tornare a Roma per il Consiglio. Se preferiamo trovarci insieme, avrei un grandissimo piacere; appena che sarà fissata l'epoca dell'adunanza di quel mese te ne scrivo: se non spero di vederti a Pisa, quando da Roma andrai in Piemonte.

Salutami Brioschi e Casorati: ricevi i saluti di Dini e di Padova ed ama sempre

il tuo Amico aff
Enrico Betti

**48**
*Carta intestata:* R. Scuola Normale Superiore di Pisa

Pisa 27/5/72 *[?]*[82]

Caro Amico,

Mi scrive il Blaserna che il 2 Giugno deve essere conferito il premio di £3000 ai Professori degli Istituti Tecnici. Della Commissione facciamo parte tu ed io. Ho veduto le pubblicazioni di 3 concorrenti: Ascoli, *[Porecini]* e Caldarera. Degli ultimi due non mi pare che sia luogo a parlarne quanto all'Ascoli, Battaglini ed io già ritenemmo meritevole d'inserirsi negli Atti la sua Memoria,[83] ed io se ne vai d'accordo non avrei difficoltà di aggiudicarli il premio. Bisogna che ti avverta che io non potrò venire a Roma il 2 Giugno; io conto che tu ci sarai, molto più che se ti occorrono schiarimenti sopra il lavoro dell'Ascoli te li potrà dare Battaglini, che ha fatto la relazione sopra il medesimo.

Tanti saluti cordiali del

Tuo Amico
Betti

---

[82] Di questo premio sembrano trattare anche le due lettere n. 68 e n. 82 che sono però datate 1878. La data di questa lettera è molto probabilmente errata, anche perché l'articolo di Ascoli viene pubblicato nel 1878: si veda la nota seguente.

[83] G. Ascoli, "Sulla rappresentabilità delle funzioni a due variabili per serie doppia trigonometrica. Memoria approvata per la stampa negli Atti dell'Accademia nella seduta del 1 giugno 1879", *Atti della R. Accademia dei Lincei*, s. III, 1879, pp. 253-300. Il verbale di tale approvazione si trova in *Atti della R. Accademia dei Lincei. Transunti*, s. III, v. 3, 1878-1879, p. 201: "L'argomento di questo lavoro è di grande importanza. L'autore intraprende per le funzioni a due variabili ricerche analoghe a quelle compiute da Riemann per le funzioni d'una sola variabile nella Memoria *Ueber die Dartstellbarkeit einer Function durch eine trigonometrische Reihe*, e giunge, in questo campo finora pochissimo esplorato, a risultamenti che ci sembrano degni d'attenzione. Perciò crediamo di poter proporre all'Accademia la stampa della presente Memoria nei suoi Atti in esteso."

**49**
*Carta intestata:* Ministero della Istruzione Pubblica - Consiglio Superiore

Roma, 10 Novembre 1872

Caro Cremona,

Bonghi mi ha detto che probabilmente la vostra Commissione sarà adunata il 17 o il 20 di questo mese.
Brioschi mi ha risposto; espone le ragioni per le quali non sia potuto far niente fin ad ora; che avrebbe intenzione di mandare una circolare colla quale sottoporrebbe alcune proposte ai membri della Società; ma prima vorrebbe esser certo che noi gli daremo appoggio. Altrimenti rinunzierà alla Presidenza. M'incarica di comunicarvi questa sua lettera. Qui non la ho; appena tornato a Pisa ti scriverò un estratto della lettera.
Il Consiglio ieri ha proposto il D'Ovidio a Professor Straordinario. Appoggerò la proposta di Armenante a Professor Straordinario, e del Bertini a incaricato di Geometria Descrittiva nella Università di Roma.
Mi dispiace che non mi troverai in Roma, ma spero di trovartici ai primi di Dicembre alle altre sedute del Consiglio; allora parleremo di molte cose.
Ama sempre

il Tuo Aff Amico
Enrico Betti

**50**
*Carta intestata:* R. Scuola Normale Superiore di Pisa

Pisa li 27 Febbrajo 1873

Mio carissimo Amico,

Non mi era dimenticato di tastare il terreno intorno a quello di che avevamo parlato a Torino. Trovai però che non so credere possibile di applicare quell'articolo della legge del 59; se fosse stato creduto possibile, sarebbe già applicato per Roma ad alcuni che non si è potuto ritenere.
Quanto alla Cattedra di Napoli, se viene l'affare al Consiglio, io sono disposto a combattere la proposta del Sannia. Cattedre di questa natura o bisogna lasciarle scoperte o darle a persone capaci e vogliose di lavorare. Quanto a Te io avrei molto più piacere, e crederei anche molto più utile che venissi a Pisa; ma se questo non è possibile non mi pare che sarà molto difficile farti chiamare a Napoli. Ora quando vo a Roma per la prossima adunanza del Consiglio, io non mancherò di occuparmene, e ti terrò informato di quelle che potrò fare. Sarò molto lieto se potrò influire nel soddisfare i tuoi desideri. [84]
La Commissione d'inchiesta[85] verrà in Toscana? vi sarà possibilità che una volta ci combiniamo in Roma?
Salutami i comuni amici della Commissione, ed ama sempre

il tuo Amico aff
Enrico Betti

---

[84] Si vedano [Menghini, 1996], lettere n. 78 e n. 79 nelle quali Cremona esplicita i motivi per i quali è scontento di restare a Milano e vorrebbe spostarsi a Napoli o a Pisa.
[85] Si tratta della Commissione d'inchiesta sull'istruzione secondaria.

**51**
*Carta intestata:* Ministero della Istruzione Pubblica - Consiglio Superiore

Roma, 10 Marzo 1873

Caro Cremona,

Oggi ho riferito sopra la proposta della Facoltà di Napoli. Nel rapporto ho detto chiaramente che una Cattedra di quella natura bisogna ricoprirla con persona che promuova quella Scienza con i suoi lavori. Il Consiglio ha approvato la mia proposta di non ammettere la dimanda della Facoltà di Napoli. Prima di procedere oltre sarebbe stato bene che parlassimo insieme; ma quando sei a Roma tu non ci sono io, quando ci sono io non ci sei tu, e perdo quasi la speranza che possiamo combinare.

Io come tu sai desidererei moltissimo di averti a Pisa, sarebbe utile per me; mi pare che sarebbe utile anche per la Scienza. Quindi vorrei fermarti quando ti sei messo in movimento a mezza strada. Almeno io non voglio aver poi a rimproverare a me stesso di non essermi adoperato a procurare a me e all'Istituto a cui appartengo un tal vantaggio. Quindi ti prego ad aprirmi l'animo tuo in questo proposito.[86]

Scrivimi a Pisa dove sarò Giovedì mattina. Ama sempre

il Tuo Amico
Enrico Betti

**52**
*Carta intestata:* R. Scuola Normale Superiore di Pisa

Pisa li 14 Settembre 1873

Amico carissimo,

Aveva saputo da Casorati del ritardo nel ricevimento del pacco che conteneva i lavori da esaminarsi per il concorso. Io avrei desiderato che negli ultimi giorni di Ottobre ci fossimo potuti riunire per trattare del giudizio da proferirsi, ma vedo che non sarà facile che per quell'epoca sia finito l'esame dei lavori dei 23 concorrenti. È vero che sono passati già vari mesi da che fu chiuso il concorso, ma un mese più o un mese meno non credo potrà essere di molto danno. Siamo soltanto cinque i membri della Commissione e tu li conosci tutti. Ho mandato i manoscritti di Geometria metà a te e metà a Bertini; quelli di Algebra a Casorati Beltrami e a me stesso.

Credo che avrai ricevuto da <u>Lemonnier</u> la prima delle tre rate degli utili del 1872 in Lire 130; la seconda la devi riscuotere i primi di Settembre. Pregherò Dini a mandare la sua Memoria; ma disgraziatamente è tornato nel Municipio ed è tutto preoccupato della risoluzione di un problema molto difficile; il rimettere in buono stato le finanze rovinate di quel Municipio.[87]

Io ho sperato più volte di combinarti a Roma, ma non è stato mai possibile, che io arrivava quando tu eri partito; l'avrei desiderato anche per parlare del progetto relativo a Napoli. Il Padula fece grandi e vivi lamenti per la decisione del Consiglio relativa alla deliberazione che prese. Io sapeva che era un uomo che non si occupava più della Scienza, ma credeva che avesse conservato il sentimento dell'interesse scientifico della sua Università.

Io sono in campagna, e spero di non essere obbligato ad andare a Roma fino ai primi di Ottobre.

A che punto siete colla vostra inchiesta? Hai pensato all'elezione del Presidente della Società Italiana?[88] Bisognerebbe valersi di questa occasione per renderle un po' di vita, ma mi pare difficile. Hai veduto le ultime riviste di Bellavitis?[89] È un gran difetto volere scrivere di ciò che non si è bene inteso.

---

[86] Si veda [Menghini, 1996], lettera n. 79, del 12 marzo.
[87] Ulisse Dini partecipò attivamente alla vita del Comune di Pisa, sia come consigliere che come assessore, per diversi anni.
[88] La carica di Brioschi è in scadenza.
[89] Serie di articoli "Rivista di giornali" che Bellavitis pubblicava sugli *Atti dell'Istituto Veneto di Scienze, Lettere ed Arti* in qualità di recensore.

Scrivimi presto, ad ama sempre

il tuo Amico aff
Enrico Betti

P.S. Mi scrisse Genocchi pochi giorni fa. È molto contento del D'Ovidio

**53**
*Carta intestata:* R. Scuola Normale Superiore di Pisa

Pisa, 28 Marzo 1874

Mio carissimo Amico,

Io sarò a Roma la mattina del 13 Aprile. Ho dato una scorsa alla traduzione del Favaro. Veramente fa vergogna a un Professore. Questo Favaro non so quando ne da chi sia stato fatto Professore, perché non sono stati esaminati dal Consiglio i suoi titoli. Del resto ne parleremo a voce. Il Meneghini[90] mi ha detto che probabilmente il Beltrami andrà Professore nella Università di Padova. È vero?

Salutami Beltrami e rammenta la Relazione sul concorso ai premi per i trattati di Matematica elementare tanto a Lui quanto a Bertini, perché quando vengo a Roma si possa finire questa faccenda.

Dini ti avrà scritto che ha respinto le prove di stampa che aveva soltanto da due giorni. Ama sempre

il tuo aff Amico
Enrico Betti

**54**
*Carta intestata:* Ministero della Istruzione Pubblica - Consiglio Superiore

Roma, 2 Settembre 74

Carissimo Amico,

Ieri il Consiglio decise intorno al Concorso di Milano; fu proposto l'Ascoli a unanimità. Seppi da Zanfi[91] l'esito dell'Affare Beltrami. Ne parleremo a voce.

Klein passò da Pisa il 25 ed io era a Pisa il 26. Me ne dispiacque molto. Io tornerò a Roma i primi di Ottobre, temo però a quanto mi dici, che per allora sarà già partito. Gli manderò una copia della mia Memoria sulla connessione degli spazi. Mi farebbe molto piacere poter parlare con Lui di questo e di molti altri soggetti.

Nell'Ottobre spero che tu sarai meno occupato e potremo stare un poco insieme.

In fretta, ti saluto caramente

il tuo aff[mo] Amico
Enrico Betti

---

[90] Giuseppe, Rettore dell'Università di Pisa dal 1871 al 1879.
[91] Si veda [Menghini, 1996], lettera n. 84, p. 66.

**55**

Sojana, 13 Settembre 74

Mio caro Amico,

Non ti ho risposto subito perché ho soltanto oggi trovato qui in campagna la tua risposta venendo da Torino.

Et sur tout pas de Zèle diceva Tayllerand [sic!] ai suoi diplomatici. Era proprio questo il caso di rammentare al Ministero questa massima. Per troppa smania di non far perdere a Te e alla Università di Roma il Beltrami non si poteva far peggio. Era già molto difficile o quasi impossibile dare una negativa recisa a un uomo del carattere e del valore di Beltrami che adduce motivi di salute e di famiglia: ma rispondergli accompagnando la negativa con una partaccia (come si dice in Toscana) era alla difficoltà della cosa aggiungerne una molto maggiore: la quistione di dignità, era obbligarlo a tutelare non solo il proprio decoro ma quello dei Professori e degli uomini che hanno dato lustro al paese con i propri lavori. Per far cosa grata a te, ti hanno dato un dispiacere.

Così non procedette il Ministero col Felici.

Riconobbe il danno che ne sarebbe venuto quando egli da Pisa fosse andato a Firenze e lo pregò caldamente a desistere dal suo progetto. Egli insistette; il Ministero mandò l'affare alla Facoltà, la quale si oppose fortemente; nonostante avrebbe acconsentito; se il Felici mosso dalla insistenza mia, della Facoltà e da ulteriori riflessioni sul suo interesse non avesse da sé stesso desistito dal suo progetto.

Io ho sempre guardato con grandissima compiacenza l'amicizia che tu e Beltrami avete l'uno per l'altro, amicizia nata e nutrita per la vita comune delle vostre belle intelligenze. È per questo che appena seppi del progetto di Beltrami ne rimasi molto male impressionato temendo che questa amicizia ne fosse turbata. Ma io non dispero che potrà sempre scansarsi un tal cattivo effetto conoscendo la lealtà di ambedue e la stima che avete reciprocamente.

La prima cosa che tu devi fare è di richiamare il Ministero sulla buona via. Beltrami è un uomo di carattere, e di una bella riputazione scientifica. È altamente vergognoso per un Governo il dimenticarlo o non saperlo. Se la risposta alla Facoltà di Padova fosse pubblicata gli tirerebbero addosso il biasimo di tutte le persone che amano la Scienza.

Quanto a me puoi star tranquillo che butterò acqua sul fuoco, quanta potrò.

Mi ha scritto Klein per combinare il modo di ritrovarsi in Toscana. Io gli scriverò a Roma.

La prima adunanza dell'Ottobre fu fissata per il 3 alle 12. Io sarò a Roma la mattina del 3 e alle 11 al Caffe del Parlamento a far colazione.

Salutami Klein, tienimi informato dell'andamento di questo affare malaugurato ed ama sempre il

Tuo Amico Aff.mo
Enrico Betti

**56**

Soiana, 21 Settembre 1874

Caro Amico,

Ho avuto molto piacere che il Ministero abbia scritto una lettera conveniente a Beltrami, ma se questa ha rimediato alle nuove difficoltà create dalla prima risposta, non ha tolto nessuna di quelle inerenti alla natura di questo affare. Io dubito molto che una parte delle difficoltà derivi dal modo differente con cui tu e Beltrami vedete la cosa. A me pare che tu sia convintissimo che Beltrami ha torto, che fa una specie di offesa a te, di tradimento alla tua Scuola, che è un uomo che si è lasciato raggirare. Egli invece che per la sua indole non può risentire i vantaggi del soggiorno di Roma, e invece ne sente gravi inconvenienti, è convinto di non fare cosa biasimevole a procurarsi un posto che crede molto più confacente alla sua indole, ai suoi studi, ai suoi interessi. Io so che ambedue questi apprezzamenti sono coscienziosi, e alla mia sincera amicizia per tutti due basterebbe questo solo, senza che m'importasse punto di giudicare del loro valore relativo. Ma per il miglior esito dell'affare io credo molto più dannoso il tuo modo di vedere del suo. Egli convinto di non aver mancato è molto più libero nelle sue decisioni; ma quali sono gli effetti che produce il tuo apprezzamento? Supponiamo anche che il Beltrami per non spingere la cosa agli estremi, sia obbligato a restare a Roma. Alle ragioni che lo hanno indotto a desiderare il trasferimento ne sarà aggiunta un altra molto più forte: la falsa posizione rispetto a un antico e caro amico che ha giudicato con severità a suo credere molto ingiusta la sua condotta nei rapporti scambievoli dell'amicizia, ed è stato cagione che venga condannato a domicilio coatto. Come Egli non dovrà desiderare molto più che per il passato di andarsene, e come a un uomo del valore di Beltrami potrà il Governo non finire per dare una tale soddisfazione?

Io ho creduto di sottoporre alla tua lucida intelligenza queste osservazioni; e come tuo amico sincero ti prego a confidarle con calma e freddezza. Avrei piacere che tu mi scrivessi il tuo parere in proposito anche prima che io venga a Roma.

Ho scritto a Klein che sarò a Firenze il 26: fammi il piacere di dirgli che non il 26, ma vi sarò il 27.

Ama sempre

il tuo aff$^{mo}$ Amico
Enrico Betti

## 57
*Carta intestata:* Ministero dell'Istruzione

Roma, 27 2/75

Caro Amico,

Il Ministro[92] dopo aver detto quelle parole alla Camera scrisse subito al Prefetto in conformità delle medesime. Quando poi seppe dell'autorizzazione che tu avevi avuto dai Ministri precedenti non mancò di scriverne al medesimo Prefetto. Essendo in appresso interpellato dal Prefetto se credesse conveniente di mantenere all'ordine del giorno la tua dimanda, Egli rispose che facesse come credeva; ma che però l'avvisava che Egli tra non molto avrebbe diretta una dimanda simile ma in una maniera determinata alla stessa Provincia; come è già stabilito nei nuovi ordinamenti degl'insegnamenti Universitari che si vanno facendo.

A me pare che in tutto ciò non vi sia nulla di personale verso di Te. La tua responsabilità è coperta da quella del Ministro, come dici benissimo. Un Ministro chiede sotto una data forma e un altro crede di dovere chiedere sotto un altra; perché Egli ha visto non solo gl'interessi di una Scuola ma quelli di tutto l'insegnamento.

D'altra parte puoi star certo che il Ministro ha alla stima del tuo valore scientifico, e che io ho per te quella sincera amicizia che ho sempre avuto da tanti anni, e la mia posizione ufficiale non può davvero mutare in nulla nessuno dei miei sentimenti. Sarebbe davvero crudele che essa oltre al sacrificio che m'impone di occupazioni tanto meno grate delle scientifiche, mi dovesse anche togliere l'amicizia degli uomini che da tanto tempo stimo ed amo. Ma io son certo che questo pericolo non vi è.

Dunque per carità considera le cose con più calma; non temere che nessuno non voglia avere riguardo alla tua dignità e conservami la tua antica amicizia.

il tuo Amico
Betti

---

[92] Ruggiero Bonghi.

**58**
*Carta intestata:* R. Scuola Normale Superiore di Pisa

Pisa 25/11/76

Caro Amico,

Io verrò a Roma Giovedì prossimo 30 del corrente mese. Avrò molto piacere di vederti e parlar teco di molte cose.

Quanto alla Memoria dal Caporali io non ne so nulla. Quando verrò a Roma ne chiederò al Segretario. Se l'avessi avuta per riferirne non avrei mancato di farlo subito.[93]

Per la seduta del 3 dei Lincei avrò da leggere una Memorietta del Roiti;[94] e forse una Nota del Dini e una mia.

Salutami gli amici ed ama sempre

il tuo Betti

**59**

Pisa 4/4/77

Caro Amico,

Ho risoluto in generale il problema che Hesse si era proposto nel V. 74 di Crelle per il caso particolare di tre corpi, e che in questo caso particolare era stato risolto da Lagrange un secolo avanti, senza la inesattezza notata da Serret nel metodo di Hesse. Dò in sostanza un metodo nuovo per l'integrazione dell'equazioni della Dinamica, che oltre a ridurre immediatamente il numero delle integrazioni al minimo, conduce ad equazioni che si prestano a trattare questi altri due problemi che ora son dietro a risolvere: determinare in quali casi un sistema di corpi che si attraggono secondo la legge di Newton si muove come un corpo rigido? in quali casi le distanze medie mantenendosi costanti, i punti si muovono periodicamente intorno alle loro posizioni medie? Intanto mi parrebbe conveniente di fare una comunicazione all'Accademia dei Lincei, riserbandomi a dare l'equazioni e le loro dimostrazioni, insieme colle applicazioni in una Memoria per gli Annali.[95]

Ti prego pertanto a far tu questa comunicazione in mio nome nella seduta di Domenica prossima, non potendo venirci a farla da me.

Le sedute del Consiglio che dovevano aver luogo il 9 sono state differite. Crederei però che il ritardo non dovesse essere di molto, e che potrò presto rivederti. Salutami i comuni amici ed ama

il tuo Amico
Betti

---

[93] Si veda [Menghini, 1996], lettera n. 89, p. 70. Ettore Caporali aveva presentato alla Società dei XL una nota per ottenere il premio annuale del valore di 400 lire, premio che gli verrà conferito nel 1878.

[94] A. Roiti, "La velocità teorica del suono e la velocità molecolare dei gas. Nota del dott. Antonio Roiti prof. nell'Istituto Tecnico di Firenze presentata dal socio E. Betti nella seduta del 3 dicembre 1876", *Atti della R. Accademia dei Lincei. Memorie*, v. I, 1877, pp. 39-45; nella stessa seduta Betti non fa altre comunicazioni. Probabilmente la Nota di Dini di cui parla in questa lettera è: U. Dini, "Sopra una classe di funzioni finite e continue che non hanno mai una derivata", che Betti presenta il 4 febbraio 1877 (*Atti della R. Accademia dei Lincei. Transunti*, v. I, 1876-77, pp. 70-72), pubblicata con il titolo "Su alcune funzioni che in tutto un intervallo non hanno mai una derivata", *Annali di Matematica*, s. II, V. VIII, 1877, pp. 121-137.

[95] E. Betti, "Sopra il moto di un sistema di un numero qualunque di punti che si attraggono o si respingono tra loro", *Annali di Matematica*, s. II, V. VIII, 1877, pp. 301-311. L'articolo di Hesse al quale fa riferimento è: O. Hesse, "Über das Problem der drei Körper", *Journal für die reine und angewandte Mathematik*, 74, 1872, pp. 97-115.

**60**
*Carta intestata:* R. Scuola Normale Superiore di Pisa

Pisa 13/5/77

Caro Cremona,

Sono passate già quasi due settimane da che mandai a Turazza i documenti del Concorso di Meccanica Razionale[96] pregandolo di esaminarli sollecitamente. Non ne ho saputo più niente. Doveva mandarli dopo al Ministero perché li rimettesse al Padula.

Sperava di vederti a Roma; ed avendomi detto il Beltrami che aveva appuntamento con te e Casorati alle 8 al Caffè del Parlamento, vi andai e aspettai sino circa alle 9; ma inutilmente. Non potei con mio dispiacere vedere né te né il Casorati. Quest'altra volta se non vengo ai Lincei bisognerà fissare in qualche modo di trovarsi insieme.

Ama sempre il

Tuo Amico aff.mo
Enrico Betti

**61**
*Carta intestata:* R. Scuola Normale Superiore di Pisa

Pisa 4/6/77

Caro Amico,

Debbo riferire al Consiglio Superiore di pubblica istruzione sopra due dimande, una dell'Ingegnere Silvio Canevazzi, l'altra del Dottore Antonio Silvani dirette la prima ad ottenere la libera docenza nelle Costruzioni, e la seconda nella Meccanica applicata. Del primo credo che tu abbia piena conoscenza essendo stato assistente nella tua scuola, del secondo forse avrai avuto occasione di conoscerlo quando eri a Bologna. Ti prego perciò a farmi il favore di scrivermi il concetto che ti sei formato del loro valore. Avrei bisogno di questa informazione prima di Domenica prossima, perché vorrei scrivere la relazione prima di partire per Roma, e credo che dovrò essere a Roma il dì 11. Anche questo mese non ho potuto intervenire ai Lincei. Quanto al concorso per la Meccanica razionale a Roma, mi ha scritto il Turazza che prima del 18 di questo mese gli sarebbe impossibile d'intervenire all'Adunanza della Commissione e mi ha pregato a differire a quella epoca questa riunione. Io penso dunque che si potrà riunirci il giorno avanti alla adunanza del Consiglio che avrà luogo i primi di Luglio e potremo in quelle sedute sbrigare la faccenda. Se il Padula non avrà sbrigato Lunedì l'esame dei Documenti, e non può sollecitamente bisognerebbe indurlo a rinunziare, e nominare subito un altro nelle prossime sedute del Consiglio.

È molto tempo che non ci siamo veduti. Questa volta conto di vederti in ogni modo e di potere discutere teco di molte cose. Ama sempre il

Tuo Amico
Enrico Betti

---

[96] Si tratta del Concorso per la nomina di un professore straordinario di Meccanica razionale nell'Università di Roma che fu vinto da Valentino Cerruti, sostenuto anche da Cremona con la lettera del 24 dicembre 1876 (si veda [Menghini, 1996], lettera n. 92, p. 72).

**62**
*Carta intestata:* R. Scuola Normale Superiore di Pisa

Pisa 29/6/77

Caro Amico,

Ho scritto subito a Turazza. Appena avuta la risposta fisserò il giorno e il luogo dell'adunanza. Quanto ad Antignano[97] io non posso altro che incitarti a seguire la buona idea che ti e venuta. Io vi fui due anni di seguito insieme col Mossotti e con i signori Bicchierai, e mi ci trovai benissimo. È un luogo di bagnatura e di villeggiatura. Gli alloggi vi sono a prezzi non esagerati. Vi sono delle villette dove si deve star assai bene. Si può passeggiare lungo il Mare e in collina secondo che piace meglio. Vi sogliono andare Caruel e Tassinari.[98] Quest'anno vi andrà la famiglia del Richiardi; egli *[però]* rimarrà molto in Pisa dove ha molto a fare per il Gabinetto e per le sue pubblicazioni. Padova sta a Livorno e verrà a trovarti; Dini a Pisa, ed io a Soiana, ma qualche corsa ad Antignano potrò farla. Aggiungi che il Ministro[99] ha nominato per il conferimento dei posti all'estero per la Matematica una Commissione di me, te e il Dini e sarà molto comodo di trovarsi vicini per il giudizio che dobbiamo proferire. Io non ho veduto ancora Bertini, ma mi dice il Padova che lo ha veduto stamattina; che rimarrà qui Luglio e Agosto e solamente 15 giorni di questi due mesi li passerà fuori. Del resto Egli stesso te ne scriverà. Vedrò di fare l'adunanza per i concorsi di Meccanica qui e allora potrai vedere da te e fissare. Quanto al giudizio per la Meccanica mi pare che saremo pienamente d'accordo. Quanto al Silvani mi pare un poco difficile di contentare il Razzaboni. Io gli aveva suggerito di nominarlo assistente e fargli fare il corso come assistente. Dopo aver dato prove nell'insegnamento sarà più facile dargli la libera docenza. Ne parleremo quando ci troveremo insieme.
Tanti saluti degli amici di qua, e del

tuo Amico
Betti

**63**

Soiana, 29/7/77

Caro Amico,

Avrai ricevuto dal Dini una lettera nella quale ti rendeva conto del colloquio avuto a Viareggio col Ministro.[100] Io aggiungerò qualche altra informazione. Parlai del modo imaginato *[sic!]* per effettuare la cosa, ed Egli mi suggerì che facessi interpellare dal Rettore[101] il Bertini se fosse disposto a passare dalla Geometria Superiore alla proiettiva, e nel caso affermativo nelle proposte che ora deve fare delle conferme degli straordinari, scrivesse al Ministero che nel caso che il Ministro potesse ricoprire la Cattedra di Geometria Superiore con un uomo eminente, il Bertini non avrebbe difficoltà ad accettare l'insegnamento di Geometria proiettiva. Io scrissi subito al Rettore, il quale mi rispose che avrebbe fatto quanto gli suggeriva molto volentieri. Parlai poi con Bertini e gli dissi che avrebbe ricevuto questa interpellanza dal Rettore, ed Egli mi promise di rispondere affermativamente; aveva già saputo che il Rettore aveva chiesto alla Segreteria dove avrebbe dovuto indirizzare una lettera a Bertini.
Quanto all'articolo 73 disse di averci già pensato da sé, perché riconosceva che tu ci scapitavi assai; però non voleva fare contemporanei i due Decreti perché non paresse che ti premiava per aver lasciato Roma. Dopo poco tempo dal tuo trasferimento manderebbe la proposta dell'applicazione dell'articolo 73, al Consiglio Superiore. Mi disse ancora che quando ricevé la tua lettera nella quale gli dicevi che avevi da parlare con Lui Egli aveva indovinato che cosa gli avresti detto. Parlò sempre come la cosa fosse decisa e non ci vedesse alcuna difficoltà. Disse però che il Decreto di trasferimento l'avrebbe fatto, insieme con quello che nominava il tuo successore alla Direzione della Scuola di Applicazione, cioè verso i primi di Novembre, perché anche per nominare Reggente il Favero gli pareva necessario che fosse Professore ordinario, ed a questo ci voleva un poco di tempo. Per ora non mi pare che ci si da fare nulla. Se però tu avessi da suggerirmi qualche cosa scrivimelo. Io ne dimandai a Coppino, ed Egli mi rispose che tu avevi già fatto quanto era necessario. Egli aveva la tua dimanda verbale e questa bastava.

---

[97] Comune in provincia di Livorno. Si veda [Menghini, 1996], lettera n. 97, p 74: Cremona stava organizzando le proprie vacanze.
[98] Potrebbe trattarsi di Théodore Caruel e di Paolo Tassinari, docenti rispettivamente di Botanica e di Chimica a Pisa.
[99] Michele Coppino.
[100] Michele Coppino.
[101] Giuseppe Meneghini.

Ho avuto lettera da Beltrami che mi rammenta l'affare Boschi. Io ho ricevuto la risposta anche da Battaglini. Mi mancano la tua e quella di Sannia. Mandai due giorni fa la relazione sui posti di perfezionamento. È bene che facciano il Decreto per Caporali, perché gli sarà di titolo, e perché potranno poi convertire quel posto in due all'interno, e noi potremo avere due bravi allievi in più. Mi pare che ciò rientri nelle idee di Coppino. A De Paolis forse non sarà ostacolo la malattia della madre a venire a Pisa.

Scrivimi ed ama sempre

il tuo Betti

**64**
*Carta intestata:* R. Scuola Normale Superiore di Pisa

Pisa, 27/9/77

Caro Cremona,

Sabato alle 8 antimeridiane verrà alla Minerva il Prof$^{re}$ Dini; alle 9 e 6 minuti partirete insieme per Pontedera e verrete a trovarmi a Soiana; di dove Domenica il giorno potrai partire per Firenze, e col treno della sera partire per Roma, ed essere la mattina del 1° ottobre alle 6 e ¼ a Roma.

Addio a Sabato

il Tuo Amico
Betti

**65**
*Carta intestata:* R. Scuola Normale Superiore di Pisa

Soiana 11/10/77

Caro Amico,

Non ti dirò quanto mi abbia fatto dolere la notizia che mi hai dato per non tormentarti inutilmente. Io capisco la difficoltà che vi era a resistere, e non ti terrò il broncio, ne diminuirà la mia sincera amicizia per te, perché non hai potuto non cedere.[102] Le adunanze del Consiglio cominceranno il 18. Dunque presto ci vedremo e potremo discorrere a voce. Vedrò Lunedì prossimo il Dini a Pisa, e gli farò leggere la tua lettera. Ama sempre il

tuo Amico
Betti

P.S. Per i posti all'interno sono due soli concorrenti, dei quali uno è il *[M...]* di cui abbiamo veduto i documenti per i posti all'estero. Potremo sbrigarcene facilmente in poco tempo.

---

[102] Si veda [Menghini, 1996], lettera n. 105, p. 79 nella quale Cremona spiega a Betti i motivi per cui ha dovuto rinunciare al trasferimento a Pisa.

**66**
*Carta intestata:* R. Scuola Normale Superiore di Pisa

Pisa 6 del 78

Caro Cremona,

Quando fui a Roma avanti Natale fui dispiacente di non poter vederti, ma non ebbi tempo di venire a trovarti come avrei desiderato. Spero di vederti quest'altra volta che sarà, se il Ministro[103] non ci chiama prima, verso la fine del mese.

Ti ringrazio della presentazione fatta in mio nome di quei teoremi di Meccanica.[104] Ho veduto nell'ultimo fascicolo dei <u>Mathematische Annalen</u> una memoria di Meyer *[sic!]* sopra lo stesso soggetto.[105] Egli determina le condizioni da verificarsi affinché sussistano i tre principj generali: io invece ne trovo dei più generali e dò le condizioni perché siano verificati. Sicché i suoi teoremi non sono altro che casi particolari dei miei, e che io stesso ho notato. Nella settimana entrante il tipografo Nistri di Pisa incomincerà la stampa della 2ª edizione della mia <u>teorica delle forze</u> etc. L'aumenterà molto; credo che verrà un libro di più di 300 pagine.[106]

Quanto a Genocchi sarò ben contento se sarà soddisfatto il suo desiderio, e se l'affare verrà al Consiglio io sosterrò la tua dimanda. La cosa è conveniente utile e giusta sotto tutti i rapporti.

Ho fatto i tuoi saluti e auguri agli amici che te li restituiscono. Accettali e partecipane alla tua famiglia da parte mia e conservami la tua amicizia

Enrico Betti

**67**
*Carta intestata:* R. Scuola Normale Superiore di Pisa

Pisa 11/5/78

Caro Amico,

Dalla tua lettera risulta che tu sei favorevole come gli altri Commissari alla promozione del Prof. Bruno.[107] Bisogna però che la lettera rimanga in atti, e che contenga esplicita la dichiarazione dell'applicabilità dell'art. 69 della legge 13 Novembre 1859.[108] Questo articolo parla di meritata fama per scritti pubblicati, scoperte, e insegnamenti dati. Gli altri si sono fondati principalmente sugli insegnamenti dati. E anche a me pare che si possa applicare a Lui come fu fatto per Boschi, l'articolo 69 per questa ragione.

Io sarò a Roma la sera del 15. Quindi tu mi potresti fare avere la lettera a Roma al Consiglio Superiore per il dì 16.

Tanti cordiali saluti del

Tuo Amico
Betti

---

[103] Francesco De Sanctis.
[104] E. Betti, "Sopra una estensione dei principi generali della dinamica", *Atti della R. Accademia dei Lincei. Transunti*, v. II, 1877-78, pp. 32-33.
[105] A. Mayer, "Über den Multiplicator eines Jacobi'schen Systems", *Mathematische Annalen*, 12, 1877, pp. 132-142.
[106] E. Betti, *Teorica delle forze newtoniane e sue applicazioni all'Elettrostatica e al Magnetismo*, Pisa, Tip. Nistri, 1879, 359 pagine. Quest'opera venne tradotta da W.F. Meyer nel 1885. Si veda la lettera n. 18 per la prima edizione.
[107] Si tratta del passaggio di Giuseppe Bruno al ruolo di ordinario di Geometria descrittiva all'Università di Torino. Si veda [Menghini, 1996], lettera n. 107, p. 81.
[108] Si tratta della legge Casati - Art. 69. Il Ministro potrà proporre al Re per la nomina, prescindendo da ogni concorso, le persone che per opere, per iscoperte, o per insegnamenti dati, saranno venute in meritata fama di singolare perizia nelle materie cui dovrebbero professare.

**68**
*Carta intestata:* R. Scuola Normale Superiore di Pisa

Pisa 28/5/78

Caro Amico,

Ti mando i lavori dei concorrenti *[Porecini]* e Caldarera: la Memoria dell'Ascoli non la ho riavuta; potrai riceverla dalla Segreteria dell'Accademia. Ti scrissi già il mio parere.[109] Ma ora debbo aggiungerti ancora un altra osservazione. Seppi da Blaserna che tra i concorrenti vi è anche il Roiti. Se fosse giudicato anch'Egli meritevole di premio, io troverei conveniente che il premio fosse diviso tra i due. A me pare che tra i due non vi sia tal differenza da anteporre l'uno all'altro. In ogni modo poi me ne rimetto a quello che voi desiderate.
  Tanti saluti cordiali del

Tuo Amico
Betti

**69**
*Carta intestata:* R. Scuola Normale Superiore di Pisa

Pisa, 5/7/78

Caro Amico,

Ho ricevuto oggi dal Ministero i documenti per il concorso ai posti di perfezionamento all'Estero. Mi piace moltissimo la tua proposta ed io scelgo Pisa, e vi sarò molto grato se invece che in città verrete a Soiana. Ne scrivo subito a Beltrami, e mi rimetto a voi per il giorno che vi sarà più conveniente.
  Quanto ai libri di testo, se tu rassegni il mandato, sarà necessario aspettare il 24 Luglio per nominare il professore, e quello farà più presto di te il lavoro? Mi pare che bisognerà rassegnarsi ad aspettare che tu abbia tempo disponibile al tuo ritorno. Mi raccomando però che tu abbia la cosa presente e ti lasci il tempo libero per questa occupazione.
  Tanti cordiali saluti del

Tuo aff.mo Amico
Enrico Betti

P.S. I concorrenti sono tre soli: Bianchi,[110] Ricci e Pittarelli

**70**

Pisa 10/7/78

Caro Amico,

Sta bene; io dimani vo a Soiana e Lunedì (15) a ore 9 antimeridiane sarò a Pisa e vi aspetto alla Scuola Normale. In questi giorni dovendo scrivermi l'indirizzo è <u>Peccioli</u> per <u>Soiana</u> (Provincia di Pisa). L'indirizzo dei telegrammi sarebbe (Bagni di Casciana per Soiana Provincia di Pisa). Ma tra l'arrivo del telegramma e l'arrivo della lettera non vi sarà molta differenza, dovendo il telegramma esser portato da un *[...]* che fa il suo comodo.
  Tanti saluti del

Tuo Amico
Betti

---

[109] Si vedano le lettere n. 48 e n. 82.
[110] Luigi.

**71**
*Carta intestata:* Ministero della Istruzione Pubblica - Consiglio Superiore

Roma 30/4/79

Caro Cremona,

Mi dispiace moltissimo di non potere venire da te alle 9; ma debbo partire. Ci vedremo quest'altra volta che sarò a Roma. Ti consegnerà il Tonelli la nota delle soscrizioni per il Chelini colla prima raccolta.[111]
Tanti complimenti alla tua signora e a tutti di tua famiglia, e una stretta di mano del

tuo Amico aff.mo
Enrico Betti

**72**
*Carta intestata:* R. Scuola Normale Superiore di Pisa

Pisa 7/12/79

Carissimo Amico,

Ti ringrazio moltissimo della tua lettera. Ho consultato il Dini, ma non posso darti subito una risposta completa alle dimande che mi dirigi.[112]
L'Editore non vorrebbe intraprendere la stampa di un'Opera come questa, spendendo Egli solo il Capitale non indifferente che è necessario. Quindi la somma che potesse dare il Ministero dovrebbe esser data per partecipare alla spesa occorrenti per la stampa. Ma non ha fatto ancora tutti i suoi conti per poter dire la somma che richiederebbe. Il dì 20 di questo mese io verrò a Roma, e probabilmente verrà anche il Dini per la Commissione dei libri di testo, si potrebbe aspettare allora a combinare la proposta che tu potresti fare quando ne fossi richiesto? Se no; avvisamene con una Cartolina e ti darò informazioni complete sulle proposte dell'Editore e sul minimo che bisogna al Dini per potere dar seguito a questa pubblicazione.
Sono dispiacentissimo di sentire che la tua Signora non è anche ristabilita. Ti prego a salutarla da parte mia ad a voler sempre bene al

Tuo Aff.mo Amico
Enrico Betti

**73**
*Carta intestata:* R. Scuola Normale Superiore di Pisa

Pisa 9/3/81

Caro Cremona,[113]

Il Ministero ha dato alla Commissione l'incarico di giudicare se al <u>Platner</u> sia applicabile l'art. 69. Le risposte degli altri Commissari sono <u>negative</u>. Quando sia tale anche la tua non vi sarebbe per questo necessità di riunire la Commissione. Ma vi è un altra quistione da risolvere. Nel Decreto di nomina del

---

[111] Cremona e Beltrami si erano fatti promotori di una sottoscrizione per erigere un monumento a Chelini all'Università di Roma, inaugurato all'inizio dell'a.a. 1879-80. L'elenco dei sottoscrittori si trova in A. Leonetti, *Pubbliche testimonianze di Stima e Affetto alla Memoria del P. Domenico Chelini delle Scuole Pie*, Bologna, 1880. Si veda [Enea-Gatto, 2009].

[112] Si veda [Menghini, 1996], lettera n. 109, p. 83 denominata "Riservatissima" da Cremona che deve rispondere al Ministro circa la richiesta di un sussidio da parte di Dini per pubblicare il proprio Corso di Analisi. Cremona, che è favorevole, chiede indicazioni a Betti.

[113] Giacomo Platner partecipò al concorso per la cattedra di Algebra complementare e Geometria analitica presso l'Università di Pavia, ma non lo vinse. Le risposte di Cremona sono: [Menghini, 1996], lettere nn. 116 e 117, pp. 85-86.

Platner è detto tra le altre ragioni, che si nomina Straordinario per essere stato dichiarato eligibile in un concorso. Quindi pare che a Lui non sia necessario applicare l'art. 69; ma invece si possa applicare quell'articolo del Regolamento secondo il quale potrebbe esser promoso ordinario dopo tre anni di effettivo esercizio scolastico a contare dal giorno della nomina, quando avesse in questo tempo acquistato titoli nuovi scientifici e didattici da meritargli una tal promozione. Gli altri Commissari hanno fatto il giudizio anche in questa ipotesi. Quando tutti i giudizi fossero unanimi anche in questo riguardo sarebbe necessario riunire la Commissione. Soltanto se non vi è unanimità io dovrò riunire la Commissione, ed è appunto per ciò che io sollecitava la tua risposta; ma non vi è nessuno inconveniente ad aspettare che tu ne abbia parlato col Beltrami il quale come collega è in caso di avere una più completa cognizione dei meriti acquistati dal Platner negli anni del suo insegnamento.

Ricevi intanto i miei più cordiali saluti, che ti prego partecipare alla tua famiglia, e credimi sempre

il Tuo Enrico Betti

**74**
*Carta intestata:* R. Scuola Normale Superiore di Pisa

Pisa 3/2/82

Caro Cremona,

Per la medaglia di oro per la Matematica che deve essere conferita dalla Società dei XL, la Commissione divide i Candidati in due categorie: Geometria ed Analisi. I candidati proposti per la prima categoria sono: D'Ovidio, Caporali, De Paolis, Veronese, Bertini e Bianchi. Per decidere sulla graduazione io ho bisogno del parere delle persone più competenti nelle materie. Perciò mi rivolgo alla tua amicizia, pregandoti a volermi comunicare confidenzialmente il tuo parere.[114]

Ringraziandoti anticipatamente ti saluto cordialmente

Tuo aff Amico
Enrico Betti

**75**
*Carta intestata:* R. Scuola Normale Superiore di Pisa

Pisa 30/11/84

Caro Amico,

Mille grazie delle congratulazioni e della offerta amichevole.[115] Era mio desiderio di essere introdotto in Senato in mezzo ai due matematici del Senato e miei antichi amici: Brioschi e Cremona.
Ho ricevuto le bozze del progetto e la relazione. Convalidata la mia nomina verrò a postar giuramento, e allora parleremo di questo progetto.

Ricevi i più cordiali saluti del

Tuo aff.mo Amico
Enrico Betti

---

[114] Si veda [Menghini, 1996], lettera n. 121, p. 87. Nel 1882 la medaglia venne assegnata a Alfredo Capelli, mentre nel 1883 a Luigi Bianchi. Caporali l'aveva già ricevuta nel 1878, D'Ovidio nel 1879 e De Paolis nel 1880 (si veda il sito dell'Accademia dei XL).
[115] Il 28 febbraio venne annunciata in Senato la nomina di Betti a Senatore e Cremona si congratula. Si veda [Menghini, 1996], lettera n. 122, p. 88.

**76**
*Carta intestata:* R. Scuola Normale Superiore di Pisa

Pisa 26 Marzo 1886

Caro Cremona,[116]

Mi scrive Genocchi che ha ricevuto una lettera di Hermite colla quale gli chiede di fare adesione insieme con i suoi amici matematici d'Italia a una testimonianza in favore della Sig$^{ra}$ Sofia Kowalewski, alla quale hanno fatto già adesione Jordan, Darboux, Appell, Poincaré, Picard e Tisserand. Genocchi ha già aderito unitamente a D'Ovidio, e prega me a mandare la mia adesione, e a scrivere a te, Beltrami, Casorati etc per sentire se acconsentite a unirvi a noi.

La ragione di questa testimonianza a favore della Sig Kowalewski è il rigetto della proposta di nominarla Socia dell'Accademia di Stokolm *[sic!]* il quale ha destato mali umori e malevolenza verso di Lei per la quale è stato anche posto in dubbio il suo valore matematico.

Genocchi avrebbe scritto da sé, se la cateratta progrediente e altri incomodi di salute non gli rendessero troppo grave scrivere lettere.

Ricevi i più cordiali saluti del

Tuo aff Amico e Collega
Enrico Betti

**77**
*Carta intestata:* R. Scuola Normale Superiore di Pisa

Pisa 5/11/86

Caro Cremona,[117]

Prima del 1859 per avere il libero esercizio della professione d'Ingegnere dopo avere riportato la Laurea in Matematica applicata era necessario di far due anni di pratica presso un Ingegnere e subire un esame a Firenze da una Commissione d'Ingegneri governativi. Dopo il regolamento Matteucci[118] fu variato soltanto il luogo dove si adunava la Commissione esaminatrice, e la sua composizione. Il luogo da Firenze fu cambiato in Pisa, e la Commissione fu composta di Professori di scienza applicata e di Ingegneri. Questo sistema durò fin al 1875; quando fu tolto a Pisa e Bologna il diritto di fare gl'Ingegneri e fu concesso invece il 1° anno di Scuola di Applicazione.

Ridolfi[119] aveva intenzione di modificare il sistema allora vigente di far gl'Ingegneri; ma non ebbe il tempo altro che di fare un Decreto con cui stabiliva qualche cosa che avrebbe condotto se fosse stato attuato a stabilire in Toscana il sistema vigente ora in Italia. Questo Decreto si trova nella Raccolta degli Atti del Governo della Toscana dal 1859 al 1861.

Questo Decreto non solo non fu attuato, ma rimase dimenticato, perché Ridolfi poco dopo lasciò il Ministero, e la istruzione pubblica di Toscana venne sotto la direzione del Ministro Italiano, a cui nessuno fece conoscere che esisteva quel Decreto.

Mi ha fatto grandissimo piacere la buona notizia che mi dai della completa riconciliazione tra te e Brioschi. Niente poteva essere per me più desiderabile, che tanta stima e amicizia ho per ambedue.

Ricevi i miei più cordiali saluti

tuo aff Amico
Enrico Betti

---

[116] Questa lettera è stata pubblicata nell'articolo: R. Tazzioli, "New perspectives on Beltrami's life and work - Considerations based on his correspondence", in: S. Coen (a cura di), *Matemathicians in Bologna (1861-1960)*, Birkhäuser, 2012, pp. 465-518. A p. 509, nota 96, è trascritta anche la lettera di Hermite a Genocchi di cui si parla qui. La lettera di Genocchi a Betti si trova invece in: N. Palladino, A.M. Mercurio, F. Palladino (a cura di), *Le corrispondenze epistolari Brioschi – Cremona e Betti – Genocchi,* Olschki, Firenze, 2009.

[117] Questa è la risposta alla domanda posta da Cremona il quale desidera sapere: "che cosa occorresse di fare dopo la detta laurea *[in scienze matematiche]*, per aver il libero esercizio della professione d'ingegnere architetto." Si veda [Menghini, 1996], lettera n. 125, p. 89.

[118] Si tratta del *Regolamento generale delle università del regno d'Italia*, del 1862.

[119] Potrebbe trattarsi di Cosimo.

**78**
*Carta intestata:* Ministero dell'Istruzione

Roma 3/12[120]

Caro Amico,

Ti rammenterai che l'ultima sera che ci vedemmo al Morteo[121] ti dissi che aveva da parlar teco di qualche cosa. Ecco di che si trattava.

Sai che Bertini per la perdita dell'incarico è rimasto qui in qualche disagio e col desiderio naturale di trovarsi meglio. Io pensando al modo di soddisfarlo, senza però che Egli me lo dimandasse, ebbi in mente che a Pisa manca con danno degli allievi della Scuola Normale l'insegnamento della Geometria Superiore; che Egli è forse quello dei tuoi allievi che più abbia profittato e che sia meglio in grado di esporre la Scienza nel vero spirito geometrico appreso da Te ed è nato in me fortissimo desiderio di mandarlo a Pisa. Ho saputo che Egli era punto alieno da acconsentirvi. Però una cosa soprattutto mi preoccupava e questa era il lasciare incerta la surrogazione qui in Roma, ed io non voleva davvero togliere uno a Roma senza prima avere da surrogarlo con uno di valore non minore. Di questo ti voleva parlare. Mi pareva però la cosa difficile e non sperava molto che tu avessi avuto alcuno da suggerirmi. Quando accidentalmente mi venne avanti l'idea che il Dino avrebbe potuto essere un eccellente successore del Bertini. Seppi che Dino ne sarebbe stato lietissimo, ed allora ritenni la cosa da farsi. Tutte le sere aspettava di vederti per parlartene, ma vedo ormai che la stagione, la lontananza e le occupazioni ci privano per molto tempo di soddisfare al desiderio di passare un ora insieme e mi sono deciso a scrivertene; consegnando la lettera al Bertini stesso, il quale ti dirà per conto tuo le ragioni per le quali accetta volentieri la cosa.

A tuo comodo bisognerà vedersi per quella cosa dell'Agricoltura e Commercio. Di molte altre cose avrei a chiacchierare teco; ma ci vedremo probabilmente all'Accademia dei Lincei.

Ama sempre

il tuo Amico
Betti

**79**
*Carta intestata:* R. Scuola Normale Superiore di Pisa

Caro Amico,

Con Beltrami abbiamo letto il rapporto del Battaglini sopra il premio di Matematica della Società dei Quaranta. Egli esamina le Memorie presentate al concorso e trova solo prendersi in considerazione per il premio quella del Fergola ed a Lui lo conferisce.[122] Beltrami ed io non avremmo alcuna difficoltà a firmare questo rapporto, se si trattasse di un vero concorso, perché senza alcun dubbio il lavoro di Fergola è di molto e quasi senza confronto superiore agli altri presentatati. Ma si può riguardare il conferimento della Medaglia della Società Italiana come il risultato di un concorso, nel senso che si debba aggiudicare a colui che abbia presentato la miglior Memoria alla Società? Mi pare che anche tu pensassi che così non debba essere. Così non fu le altre volte, e non risulta neppure dalle parole dello Statuto. Prima però di scrivere al Battaglini perché rifaccia il rapporto voleva che ben si stabilissero i criteri con i quali il premio si conferisce, ossia voleva che tu mi dicessi se crederesti anche tu che si potesse fare il rapporto nel modo che vo a dire.

Il lavoro del Fergola non è senza importanza ma la importanza del soggetto, e il valore del lavoro non so se arriverebbe a quel grado che si richiederebbe quando si dovesse andare a scegliere il miglior lavoro da premiarsi in tutto il campo della Matematica pura ed applicata in Italia. Quindi non potendo limitare il campo alle Memorie presentate converrà forse limitarlo in altro modo. Con Beltrami si pensava se potrebbe la Commissione dire nel suo rapporto, che per ben tre volte avendo aggiudicato la Medaglia per lavori di Matematica pura, ha creduto conveniente di prendere in considerazione in questo anno i lavori di Matematica applicata, e così trovarsi più nel vero nel giudizio che Ella farà del merito della Memoria. Così dare Facoltà alle Commissioni di limitare la loro ricerca dei lavori più meritevoli un anno a una parte uno all'altra quando lo trovino utile, si renderà anche meno difficile il lavoro dell'altra Commissione che deve scegliere il migliore lavoro in un campo tanto più vasto e più eterogeneo. Scrivimi che cosa ne pensi, ed io subito scriverò subito al Battaglini, perché rifaccia il rapporto con questo nuovo punto di vista.

---

[120] Potrebbe essere del 1874: Cremona è a Roma e Bertini si sposterà a Pisa nel 1875.
[121] Il *Caffè Morteo* era un locale di Roma, in via del Corso. Fu demolito nel 1881 per far spazio alla *Rinascente*.
[122] Questo fatto e la risposta di Cremona (si veda [Menghini, 1996], lettera n. 93, p. 72) portano a datare la lettera al dicembre 1875, perché nel 1876 a Emanuele Fergola venne conferita la medaglia della Società dei XL.

    Cantoni[123] ebbe la partecipazione prima di partire da Roma. Pacinotti ha avuto anch'esso la partecipazione; ma Villari non ha ancora inviato i documenti di Pisati. Ho incaricato Beltrami che dimani lo vedrà a Bologna di sollecitarlo.

    Ama sempre il

<div style="text-align:right">Tuo Amico<br>Betti</div>

---

[123] Potrebbe trattarsi di Giovanni Cantoni.

*Lettere a Eugenio Beltrami*

**80**

Chiarissimo Sig.r Professore,[124]

Vorrebbe Ella essere nominato Professore ordinario di Geodesia teoretica sulla R. Università di Pisa? Mi scriva subito sì o no.

Il desiderio di avere un Professore che possa fare lezione di questa scienza come si conviene in una Università mi ha fatto ardito a rivolgermi a Lei sebbene non la conosca personalmente, ma solo abbia notizia della sua attitudine scientifica per i lavori che ha pubblicato nostri Annali e per quello che mi hanno detto di Lei i miei cari amici Brioschi e Cremona. Colle sue cognizioni di analisi, e con due o tre mesi che Ella stia in un Osservatorio o a Milano o a Firenze, si può porre in grado di far benissimo questa lezione; poi potrà (coltivando questo ramo di Scienza) far molto. Io no ho già parlato al Ministro[125] il quale attende una mia lettera per fare l'opportuno decreto di nomina.

Aspetto dunque sollecita sua risposta che spero affermativa e mi dichiaro

Suo aff Collega
Enrico Betti

Spezia, 11 Agosto 1861 *[1863]*[126]

**81**
*Carta intestata:* R. Scuola Normale Superiore Pisa

Pisa, li 30 Maggio 1872

Caro Cremona,[127]

Io non mi rammento le precise parole della risoluzione del Consiglio e della Relazione.[128] Mi rammento bene però che non solo furono contrarie alle dimande del petente, ma ancora che fu risposto negativamente alla dimanda del Ministero, se convenisse impiegarlo come Professore Straordinario in qualche cattedra Universitaria: e che la ragione di questa risposta fu la seguente. Egli nel 1865 fu chiamato come Professore Straordinario a insegnare Meccanica nella Università di Bologna, e il Ministero stesso ci faceva sapere che fu obbligato a non confermarlo nell'anno successivo perché non aveva soddisfatto nel suo insegnamento: che non si avevano d'allora in poi altre prove che inducessero a modificare il giudizio non favorevole che si era formato sulla sua abilità scientifica.

Questo giudizio del Consiglio superiore che io ricordando bene quello che mi diceva di Lui l'amico Cremona nel 1865,[129] credo conforme al vero, si riferisce solo all'insegnamento Universitario. Forse il pover uomo potrebbe fare discretamente bene in un Istituto tecnico ed in un Liceo ed io ben volentieri, quando anche tu che lo conosci meglio fossi di questa opinione, in quel che posso gli gioverei perché gli fosse data una tal posizione. Mi pare che se è possibile questa sia la miglior soluzione.

---

[124] Questa lettera è indirizzata a Eugenio Beltrami (si veda qui sotto la nota 126).
[125] Michele Amari.
[126] L'anno scritto da Betti è molto probabilmente errato. Si veda infatti [Menghini, 1996], lettera n. 14, p. 16 e la 049-09786 scritta il 16 agosto 1863 da Beltrami a Cremona, che inizia così: "Di ritorno questa mattina dalla campagna ho trovato la lettera del Betti che qui unita ti trasmetto." In questo modo si spiega anche il motivo per cui essa si trovi in questo Archivio. Inoltre si veda la risposta di Beltrami a Betti in [Giacardi-Tazzioli, 2012], p. 63. Beltrami accetterà la proposta e andrà a Pisa sulla cattedra di Geodesia teorica dal 1863 al 1866.
[127] Questa lettera non è indirizzata a Cremona, ma a Beltrami: si veda infatti la risposta di Beltrami a Betti del 23 giugno in [Giacardi-Tazzioli, 2012], p. 98. Inoltre Cremona viene citato più volte nel testo e si trovava a Milano, mentre Magni era a Bologna dove era anche Beltrami che avrebbe potuto così presentargli i saluti di Betti.
[128] Si parla di Luigi Venturi.
[129] Si vedano [Menghini, 1996], lettere nn.36 e 37, pp. 28-29.

Quanto al Concorso appena avrai finito l'esame dei Documenti, potrai passarli a Cremona perché ne termini l'esame se avrà terminato la sua ispezione; se poi, come credo, finirai di esaminarli prima rimettili a Tardy, pregandolo a passarli a Cremona dopo averli esaminati. Quanto alle dissertazioni dei concorrenti per esame basta che le presentino in tempo perché possa essere distribuita almeno 8 giorni prima dell'esperimento orale ai Commissari e mi pare che ci sia anche tempo.

Io avrei avuto molto piacere di poterne uscirne nel Giugno, ma se non sarà possibile si farà nel Luglio dopo terminati gli esami, o quando vi farà più comodo. Disturberà un poco, ma come si fa?

Ricevi i saluti di Dini Felici e Padova. Salutami Cecco Magni ed ama sempre

<div style="text-align:right">il tuo Amico<br>Enrico Betti</div>

*Lettera a Pietro Blaserna*

**82**
*Carta intestata:* R. Scuola Normale Superiore di Pisa

Pisa 28/5/78

Caro Blaserna, [130]

    Ho scritto al Cremona; non ho ancora avuto risposta; ma a quanto mi scrisse prima di partire dovrebbe essere di ritorno l'ultimo del mese, e quindi potrà trovarsi presente all'adunanza della Commissione.
    Io trovo che si può conferire il premio all'Ascoli, e se poi la Commissione credesse che ci fosse un altro meritevole del premio, non avrei difficoltà alcuna ad acconsentire che il premio fosse diviso. Per esempio, tra il Roiti e l'Ascoli io non credo che vi sia tal differenza che non si possa dividere il premio tra i due. Anzi a me parrebbe la cosa più conveniente.
    Tanti saluti cordiali del

Tuo aff Amico
Enrico Betti

Non so se Felici potrà venire a Roma perché è infreddato.

Aveva già scritto la lettera quando ho ricevuto il tuo Dispaccio ed ho mandato subito a Cremona i due libri dei concorrenti che aveva ricevuto *[sic!]*. Ti ringrazio molto che ti sia prestato a risparmiarmi la gita a Roma che mi avrebbe molto scomodato.
    Tanti saluti agli amici di Morteo.[131]

Betti

---

[130] Si veda la nota 82.
[131] Si veda la nota 121.

*Lettera di Giorgini[132] a Betti*

**83**

Caro Betti

Firenze 1.3.76

Ti mando tutto quello che ho potuto raccogliere del povero Babbo.[133] Non ho trovato per la casa copia di un lavoro che ho sentito pure rammentare. La teoria analitica delle projezioni.[134] Unisco una biografia scritta da un Giovanni Sforza,[135] che non diffusi, non avendo il nome dell'autore nessuna autorità, e perché stonature che mi avevano offeso l'orecchio, ma nella quale si trovano *[penso]* notizie quante bastano ai fini del Cremona.[136] Vorrei riavere la Statica[137] e il volumetto del Giornale della Scuola Politecnica[138] non avendo in casa altre copie. Tante cose al Bonghi.

Il tuo Aff.
Giorgini

---

[132] Non si tratta di Giovanni Battista, perché la sua grafia è diversa.
[133] Gaetano Giorgini, morto il 16 settembre 1874.
[134] Pubblicato a Lucca nel 1820.
[135] G. Sforza, *Nelle esequie solenni del senatore G.G. celebrate nella chiesa parrocchiale di Montignoso il 23 settembre 1874. Discorso*, Lucca, 1875.
[136] In [Menghini, 1996], lettera n. 89, pp. 70-71 Cremona aveva scritto di stare terminando le notizie sul Giorgini, ma non risulta alcun necrologio firmato da lui.
[137] *Elementi di statica*, Firenze, 1835.
[138] Gaetano Giorgini studiò a Parigi presso l'École polytechnique. Probabilmente si tratta del suo primo contributo scientifico: *Démonstration de quelques théorèmes de géométrie*, nella *Correspondance sur l'École impériale polytechnique*, III, 1816, pp. 6-9.

## Tabella - dati delle lettere

|    | **DESTINATARIO** | **LUOGO E DATA** | **SEGNATURA FILE E CARTACEA** | **CONSISTENZA** (n. di pagine) |
|----|------------------|------------------|-------------------------------|--------------------------------|
| 1  | Luigi Cremona | Pisa 19 novembre 1860 | 089-20923 (17590) | 1 |
| 2  | Luigi Cremona | Torino 15 febbraio 1861 | 089-20924 (17591) | 2 |
| 3  | Luigi Cremona | Torino 12 marzo 1863 | 089-20896 (17563) | 2 |
| 4  | Luigi Cremona | Pisa 9 aprile 1863 | 089-20897 (17564) | 1 |
| 5  | Luigi Cremona | Pisa 15 aprile 1863 | 089-20898 (17565) | 1 |
| 6  | Luigi Cremona | Spezia 24 agosto 1863 | 089-20899 (17566) | 1 |
| 7  | Luigi Cremona | Spezia 29 agosto 1863 | 089-20900 (17567) | 1 |
| 8  | Luigi Cremona | Torino 4 dicembre 1863 | 089-20926 (17593) | 3 |
| 9  | Luigi Cremona | Torino 7 marzo 1864 | 089-20901 (17568) | 1 |
| 10 | Luigi Cremona | Torino 18 luglio 1864 | 089-20902 (17569) | 2 |
| 11 | Luigi Cremona | Pisa 19 maggio 1865 | 089-20903 (17570) | 2 |
| 12 | Luigi Cremona | Di campagna presso Pistoia 14 settembre 1865 | 050-10496 (7183) | 2 |
| 13 | Luigi Cremona | Pisa 8 aprile 1866 | 089-20904 (17571) | 1 |
| 14 | Luigi Cremona | Pisa 27 giugno 1866 | 050-10497 (7184) | 1 |
| 15 | Luigi Cremona | Firenze 17 agosto 1866 | 050-10498 (7185) | 2 |
| 16 | Luigi Cremona | Firenze 19 agosto 1866 | 050-10499 (7186) | 2 |
| 17 | Luigi Cremona | Firenze 29 agosto 1866 | 050-10500 (7187) | 2 |
| 18 | Luigi Cremona | Caloria 15 ottobre 1866 | 050-10501 (7188) | 2 |
| 19 | Luigi Cremona | Pisa 8 gennaio 1867 | 050-10502 (7189) | 3 |
| 20 | Luigi Cremona | Pisa 25 marzo 1867 | 050-10506 (7193) | 2 |
| 21 | Luigi Cremona | Pistoia 7 settembre 1867 | 050-10503 (7190) | 1 |
| 22 | Luigi Cremona | Pistoia 24 settembre 1867 | 050-10514 (7201) | 2 |
| 23 | Luigi Cremona | Pistoia 30 settembre 1867 | 050-10515 (7202) | 2 |
| 24 | Luigi Cremona | Pistoia 24 ottobre 1867 | 050-10516 (7203) | 1 |
| 25 | Luigi Cremona | Pisa 21 novembre 1867 | 050-10517 (7204) | 3 |
| 26 | Luigi Cremona | Pisa 11 dicembre 1867 | 050-10518 (7205) | 1 |
| 27 | Luigi Cremona | Pisa 27 dicembre 1867 | 050-10519 (7206) | 2 |
| 28 | Luigi Cremona | Pisa 14 febbraio 1868 | 050-10504 (7191) | 2 |
| 29 | Luigi Cremona | Pisa 7 marzo 1868 | 050-10505 (7192) | 1 |
| 30 | Luigi Cremona | Pisa 13 aprile 1868 | 050-10507 (7194) | 1 |
| 31 | Luigi Cremona | Pisa 2 maggio 1868 | 050-10508 (7195) | 1 |
| 32 | Luigi Cremona | Pisa 18 luglio 1868 | 050-10509 (7196) | 2 |
| 33 | Luigi Cremona | Pisa 23 luglio 1868 | 050-10510 (7197) | 2 |
| 34 | Luigi Cremona | Pisa 22 agosto 1868 | 050-10511 (7198) | 2 |
| 35 | Luigi Cremona | Pistoia 18 settembre 1868 | 050-10512 (7199) | 3 |
| 36 | Luigi Cremona | Pistoia 2 ottobre 1868 | 050-10513 (7200) | 3 |
| 37 | Luigi Cremona | Pisa 17 novembre 1869 | 050-10520 (7207) | 1 |
| 38 | Luigi Cremona | Pisa 30 aprile 1870 | 050-10521 (7208) | 1 |
| 39 | Luigi Cremona | Ponsacco presso Pontedera 18 Agosto 1870 | 050-10522 (7209) | 1 |
| 40 | Luigi Cremona | Pisa 3 novembre 1870 | 050-10523 (7210) | 1 |
| 41 | Luigi Cremona | Pisa 25 dicembre 1870 | 050-10524 (7211) | 1 |
| 42 | Luigi Cremona | Pisa 26 marzo 1871 | 050-10525 (7212) | 2 |
| 43 | Luigi Cremona | Pisa 16 giugno 1871 | 050-10526 (7213) | 1 |
| 44 | Luigi Cremona | Soiana 23 agosto 1871 | 050-10527 (7214) | 3 |
| 45 | Luigi Cremona | Pisa 12 novembre 1871 | 050-10528 (7215) | 1 |

| 46 | Luigi Cremona | Pisa 4 dicembre 1871 | 050-10529 (7216) | 2 |
|---|---|---|---|---|
| 47 | Luigi Cremona | Pisa 5 maggio 1872 | 050-10530 (7217) | 2 |
| 48 | Luigi Cremona | Pisa 27 maggio 1872 | 089-20905 (17572) | 1 |
| 49 | Luigi Cremona | Roma 10 novembre 1872 | 050-10532 (7219) | 3 |
| 50 | Luigi Cremona | Pisa 27 febbraio 1873 | 089-20906 (17573) | 2 |
| 51 | Luigi Cremona | Roma 10 marzo 1873 | 050-10533 (7220) | 2 |
| 52 | Luigi Cremona | Pisa 14 settembre 1873 | 089-20907 (17574) | 3 |
| 53 | Luigi Cremona | Pisa 28 marzo 1874 | 089-20908 (17575) | 1 |
| 54 | Luigi Cremona | Roma 2 settembre 1874 | 089-20909 (17576) | 2 |
| 55 | Luigi Cremona | Soiana 13 settembre 1874 | 089-20910 (17577) | 3 |
| 56 | Luigi Cremona | Soiana 21 settembre 1874 | 089-20911 (17578) | 3 |
| 57 | Luigi Cremona | Roma 27 febbraio 1875 | 089-20912 (17579) | 4 |
| 58 | Luigi Cremona | Pisa 25 novembre 1876 | 089-20913 (17580) | 1 |
| 59 | Luigi Cremona | Pisa 4 aprile 1877 | 089-20914 (17581) | 1 |
| 60 | Luigi Cremona | Pisa 13 maggio 1877 | 050-10535 (7222) | 1 |
| 61 | Luigi Cremona | Pisa 4 giugno 1877 | 089-20916 (17583) | 2 |
| 62 | Luigi Cremona | Pisa 29 giugno 1877 | 089-20917 (17584) | 2 |
| 63 | Luigi Cremona | Soiana 29 luglio 1877 | 089-20918 (17585) | 3 |
| 64 | Luigi Cremona | Pisa 27 settembre 1877 | 050-10536 (7223) | 1 |
| 65 | Luigi Cremona | Soiana 11 ottobre 1877 | 089-20915 (17582) | 1 |
| 66 | Luigi Cremona | Pisa 6 gennaio 1878 | 050-10537 (7224) | 2 |
| 67 | Luigi Cremona | Pisa 11 maggio 1878 | 089-20919 (17586) | 2 |
| 68 | Luigi Cremona | Pisa 28 maggio 1878 | 089-20920 (17587) | 1 |
| 69 | Luigi Cremona | Pisa 5 luglio 1878 | 089-20921 (17588) | 1 |
| 70 | Luigi Cremona | Pisa 10 luglio 1878 | 089-20921bis (17588) | 1 |
| 71 | Luigi Cremona | Roma 30 aprile 1879 | 050-10539 (7226) | 1 |
| 72 | Luigi Cremona | Pisa 7 dicembre 1879 | 089-20922 (17589) | 2 |
| 73 | Luigi Cremona | Pisa 9 marzo 1881 | 050-10534 (7221) | 2 |
| 74 | Luigi Cremona | Pisa 3 febbraio 1882 | 050-10540 (7227) | 2 |
| 75 | Luigi Cremona | Pisa 30 novembre 1884 | 089-20927 (17594) | 1 |
| 76 | Luigi Cremona | Pisa 26 marzo 1886 | 050-10542 (7229) | 1 |
| 77 | Luigi Cremona | Pisa 5 novembre 1886 | 050-10538 (7225) | 2 |
| 78 | Luigi Cremona | Roma 3 dicembre *[1874]* | 089-20929 (17596) | 3 |
| 79 | Luigi Cremona | s.l. s.d. *[dicembre 1875]* | 089-20928 (17595) | 3 |
| 80 | Eugenio Beltrami | Spezia 11 agosto 1861 *[1863]*\* | 089-20925 (17591) | 2 |
| 81 | Eugenio Beltrami\*\* | Pisa 30 maggio 1872 | 050-10531 (7218) | 3 |
| 82 | Pietro Blaserna | Pisa 28 maggio 1878 | 050-10543 (7230) | 2 |
| 83 | Enrico Betti da Giorgini | Firenze 1 marzo 1876 | 089-20729 | 2 |

\* L'anno indicato da Betti è 1861, ma si tratta certamente di un errore.
\*\* Betti inizia la lettera con "Caro Cremona", ma si tratta di un evidente errore.

## Indice dei nomi citati nelle lettere*

| Amari | Michele (Ministro) | 3, 7, 8, 9, 10, 11, 80 |
|---|---|---|
| Appell | Paul Émile | 76 |
| Armenante | Angelo | 49 |
| Ascoli | Giulio | 48, 54, 68, 82 |
| Baltzer | Richard | 32, 33, 35, 36 |
| Barberis | Giuseppe | 41 |
| Battaglini | Giuseppe | 3, 44, 45, 48, 63, 78, 79 |
| Bellavitis | Giusto | 52 |
| Beltrami | Eugenio | 6, 15, 16, 47, 52, 53, 55, 56, 60, 63, 73, 79 |
| Berti | Domenico (Ministro) | 13, 16, 17, 19, 20 |
| Bertini | Eugenio | 25, 37, 49, 52, 53, 62, 63, 74, 78 |
| Bertrand | Joseph | 32 |
| Bianchi | Luigi | 74, 69 |
| Bianchi | Nicomede | 12 |
| Bicchierai | (signori) | 4, 62 |
| Blaserna | Pietro | 48, 68 |
| Bonghi | Ruggiero | 8, 12, 49, 57, 83 |
| Boschi | Pietro | 16, 63, 67 |
| Brioschi | Francesco | 1, 11, 12, 13, 14, 18, 19, 20, 21, 22, 23, 25, 27, 28, 30, 32, 33, 34, 35, 37, 38, 39, 43, 45, 47, 49, 75, 77, 80 |
| Bruno | Giuseppe | 67 |
| Caldarera | | 48, 68 |
| Canevazzi | Silvio | 61 |
| Cannizzaro | Stanislao | 13 |
| Cantoni | Giovanni | 79 |
| Capellini | Giovanni | 5 |
| Caporali | Ettore | 58, 63, 74 |
| Caruel | Théodore | 62 |
| Casorati | Felice | 1, 19, 20, 25, 27, 32, 37, 38, 39, 40, 41, 42, 43, 45, 47, 52, 60, 76 |
| Chelini | Domenico | 2, 8, 71 |
| Clebsch | Rudolf Friedrich Alfred | 10, 34 |
| Codazzi | Delfino | 6 |
| Colenso | John William | 23, 24, 25, 29 |
| Coppino | Michele (Ministro) | 62, 63 |
| Correnti | Cesare (Ministro) | 20, 44 |
| Crelle | August Leopold | 19, 59 |
| Cremona | Luigi | 4, 35, 75, 80, 81, 82, 83 |
| Darboux | Jean Gaston | 76 |
| D'Arcais | Francesco | 42 |
| De Gasparis | Annibale | 3 |
| De Paolis | Riccardo | 63, 74 |
| De Sanctis | Francesco (Ministro) | 66 |
| Della Rocca | Gino | 26 |
| Dino (de) | Nicola Salvatore | 78 |
| Dini | Ulisse | 13, 16, 20, 22, 23, 25, 26, 29, 32, 33, 37, 40, 41, 43, 46, 47, 52, 53, 58, 62, 63, 64, 65, 72, 81 |
| Donati | Giovanni Battista | 7, 13, 20 |
| Piria | Raffaele | 8 |
| Pisati | Giuseppe | 79 |
| Pittarelli | Giulio | 69 |
| Platner | Giacomo | 27, 33, 34, 35, 73 |
| Poincaré | Jules Henri | 76 |
| Porecini | | 48, 68 |
| Razzaboni | Cesare | 62 |
| Richiardi | Sebastiano | 62 |
| Ricci Curbastro | Gregorio | 69 |
| Ridolfi | Cosimo | 2, 77 |
| Riemann | Georg Friederich Bernhard | 27, 44 |
| Roiti | Antonio | 58, 68, 82 |
| Rubini | Raffaele | 44 |
| Sacchi | | 32 |
| Sannia | Achille | 50, 63 |
| Salvagnoli | Vincenzo | 4 |
| Serret | Joseph Alfred | 59 |
| Sforza | Giovanni | 83 |
| Silvani | Antonio | 61, 62 |
| Simson | Robert | 27 |
| Sismonda | Angelo | 32 |
| Stoppani | Antonio | 32 |
| Sismonda | Eugenio | 32 |
| Sylvester | James Joseph | 11 |
| Talleyrand-Périgord (de) | Charles-Maurice | 55 |
| Tardy | Placido | 10, 81 |
| Tassinari | Paolo | 62 |
| Teza | Emilio | 20, 32 |
| Tisserand | Félix | 76 |
| Tonelli | Alberto | 71 |
| Tortolini | Barnaba | 19, |
| Turazza | Domenico | 60, 61, 62 |
| Venturi | Luigi | 17 |
| Veronese | Giuseppe | 74 |
| Viviani | Vincenzo | 24 |
| Villari | Emilio | 12, 79 |
| Volpicelli | Paolo | 38 |
| Zanfi | Luigi | 54 |
| Zinna | Alfonso | 26 |

**\*** Brevi biografie di quasi tutti i citati si trovano nel volume [Menghini, 1996] in nota alla corrispondenza Cremona-Betti, nel file *Nomi citati nelle lettere pubblicate nella collana "Materiali per costruzione delle biografie di matematici italiani dopo l'Unità"* e in altre pubblicazioni di carteggi su questo sito.
\*\* Si veda [Roero, 2013], pp. 315-318.

## Referenze bibliografiche

U. Bottazzini, "Francesco Brioschi e la cultura scientifica nell'Italia post-unitaria", *BUMI*, v. 1°, n. 1, 1998, pp. 59-78 (consultabile in rete)

C. Cerroni, L. Martini (a cura di), *Il carteggio Betti-Tardy (1850-1891)*, Mimesis Ed., Milano, 2009

[Cremona, *Opere*] = *Opere matematiche di Luigi Cremona*, 3 volumi, Hoepli, 1914-17 (consultabili in rete)

[Enea-Gatto, 2009] = M.R. Enea, R. Gatto, *Le carte di Domenico Chelini dell'Archivio generale delle Scuole Pie e la corrispondenza Chelini-Cremona (1863-1878)*, Mimesis Ed., Milano, 2009

[Giacardi-Tazzioli, 2012] = L. Giacardi, R. Tazzioli (a cura di), *Le lettere di Eugenio Beltrami a Betti, Tardy e Gherardi*, Mimesis Ed., Milano, 2012

[Menghini, 1996] = M. Menghini (a cura di), *Per l'Archivio della Corrispondenza dei Matematici Italiani - La corrispondenza di Luigi Cremona*, v. III, Università Bocconi, Milano, 1996

P. Riccardi, *Saggio di una bibliografia euclidea*, Bologna, 1887 (consultabile in rete)

[Roero, 2013] = C.S. Roero (a cura di), *Dall'Università di Torino all'Italia unita. Contributi dei docenti al Risorgimento e all'Unità*, Università degli Studi di Torino, 2013 (consultabile in rete)